\documentclass[final,square,comma,3p,12pt]{elsarticle}
\usepackage[english]{babel}
\usepackage[utf8]{inputenc}
\usepackage{amsmath,amssymb,color,geometry,graphics,subfig}
\usepackage[pdftex,bookmarks]{hyperref}

\providecommand{\graphicPath}{.}

\providecommand{\Title}{Anisotropic Mesh Adaptation for Variational
   Problems\\ Using Error Estimation Based on Hierarchical Bases}
\providecommand{\Journal}{Canad. Appl. Math. Quart. (special issue)}

\hypersetup{colorlinks=false, pdftitle={\Title}, pdfsubject={\Journal},
   pdfauthor={Weizhang Huang, Lennard Kamenski, and Xianping Li}, 
}

\providecommand{\boldx}{\boldsymbol{x}}
\providecommand{\R}{\mathbb{R}}
\providecommand{\cT}{\mathcal{T}}

\providecommand{\abs}[1]{\left|{#1}\right|}
\providecommand{\norm}[1]{\left\|{#1}\right\|}
\providecommand{\snorm}[1]{\left\langle{#1}\right\rangle}
\newcommand{\InkRed}[1]{#1}

\DeclareMathOperator{\tr}{tr}

\begin{document}
%********************************************************************
\begin{frontmatter}
%********************************************************************

%--- title, journal, keywords ---------------------------------------
\title{\Title}
\journal{\Journal}
\begin{keyword} 
   mesh adaptation \sep anisotropic mesh \sep finite element
   \sep a~posteriori estimator \sep hierarchical basis 
   \sep variational problem
   \MSC 65N50 \sep 65N30 \sep 65N15
\end{keyword}

%--- authors / address  ---------------------------------------------
\author[addressKU]{Weizhang~Huang}
\ead{huang@math.ku.edu}

\author[addressKU]{Lennard~Kamenski}
\ead{lkamenski@math.ku.edu}

\author[addressKU]{Xianping~Li}
\ead{lxp@math.ku.edu}

\address[addressKU]{Department of Mathematics, The University of Kansas, 
   Lawrence, KS~66045, USA}

%--- abstract ------------------------------------------------------
\begin{abstract}
Anisotropic mesh adaptation has been successfully applied to the numerical solution of partial differential equations but little considered for variational problems. 
In this paper, we investigate the use of a global hierarchical basis error estimator for the development of an anisotropic metric tensor needed for the adaptive finite element solution of variational problems.
The new metric tensor is completely a~posteriori and based on residual, edge jumps and the hierarchical basis error estimator.
Numerical results show that it performs comparable with existing metric tensors based on Hessian recovery.
A few sweeps of the symmetric Gau{\ss}-Seidel iteration for solving the global error problem prove sufficient to provide directional information necessary for successful mesh adaptation. 
\end{abstract}
\end{frontmatter}

%**************************************************************************
\section{Introduction}
%**************************************************************************

Variational problems arise from many areas of science and engineering and have a natural formulation with which the governing equation can be derived through minimization.
Often, a variational problem can be transformed into a boundary value problem of partial differential equations (PDEs) and solved by methods specially designed for PDEs.
Unfortunately, these methods are not designed to take structural advantage of variational problems and many researches argue that the variational formulation should be used as a natural optimality criterion for mesh adaptation; e.g., see \cite{BR96,BR01,BMM97}.
In this article, we consider mesh adaptation based on the variational formulation.

Most of the existing work for the adaptive numerical solution of variational problems employs isotropic mesh adaptation where the size of mesh elements varies from place to place according to a certain local error estimate while their shape is kept close to \InkRed{equilateral}.
On the other hand, the ability to allow the size, shape and orientation of mesh elements to vary can significantly improve the accuracy of the solution and enhance the computational efficiency; this is especially true for problems exhibiting distinct anisotropic features.
While it has attracted considerable attention from many researchers and been successfully applied to the numerical solution of PDEs, anisotropic mesh adaptation has rarely been employed for variational problems, especially when combined with a~posteriori error estimates.

The objective of this article is to investigate the application of a hierarchical basis a~posteriori error estimator  (HBEE) to anisotropic mesh adaptation for variational problems.
We use the \emph{$M$-uniform mesh} approach \cite{Huang05a} where an adaptive mesh is generated as a uniform mesh in the metric specified by a symmetric and strictly positive definite tensor, $M$.
The key to the success of the approach is to define a proper metric tensor.
Recently, Huang and Li \cite{HuaLi10} developed a metric tensor for the adaptive finite element solution of variational problems.
In the anisotropic case, $M$ is semi-a~posteriori: it involves residual and edge jumps, both dependent on the computed solution, and the Hessian of the exact solution, which provides necessary directional information for anisotropic mesh adaptation and is typically approximated using a Hessian recovery technique \cite{HuaLi10}.
On the other hand, an HBEE was used by Huang et al. \cite{HuKaLa10} to develop a metric tensor for the adaptive finite element solution of the boundary value problem of a second-order elliptic differential equation.
In this paper we combine the approaches in \cite{HuKaLa10,HuaLi10} and develop a metric tensor for variational problems based on an HBEE and the underlying variational formulation.
The obtained metric tensor is a~posteriori in the sense that it is based solely on residual, edge jumps, and an a~posteriori error estimate.
This is in contrast to most previous work where $M$ depends on the Hessian of the exact solution and is semi-a~posteriori or completely a~priori; e.g., see \cite{BoGHLS97,BoGeMo97,CaHeMP97,Frey08,Huang05a,HuaLi10}. 
Numerical results show that the new metric tensor works well for all numerical examples and is comparable in performance with those based on Hessian recovery.

The outline of the paper is as follows.
Section~\ref{sec:FE-Approximation} is devoted to the description of a general variational problem, its finite element approximation, and the general procedure of the $M$-uniform mesh approach for mesh adaptation.
In Section~\ref{sec:HBEE}, a bound is derived for the variation of the underlying functional.
It involves residual, edge jumps, and an HBEE and is minimized on $M$-uniform meshes for an optimal $M$.
Numerical examples are given in Section~\ref{sec:Numerical-Results}.
Finally, concluding remarks are given in Section~\ref{sect:Conclusion}.

%**************************************************************************
\section{Finite element approximation and mesh adaptation}
%**************************************************************************
\label{sec:FE-Approximation}

%--- Variational problem --------------------------------------------------
Consider a general functional of the form 
\begin{align}
   I[v] = \int_{\Omega} F(\boldx,v,\nabla v) d\boldx ,
      \quad \forall v \in U_0
   \label{functional-1}
\end{align}
where $F(\cdot, \cdot, \cdot)$ is a given smooth function, $\Omega \subset \R^{d}$ ($d = 1,2,3$) is the physical domain and $U_0$ is a properly selected set of functions satisfying the \InkRed{homogeneous} Dirichlet boundary condition 
\begin{align*}
  v(\boldx) = \InkRed{0} \quad \forall \boldx \in \partial \Omega.
   % \label{bc-1}
\end{align*}
In this paper, we consider only \InkRed{the homogeneous Dirichlet boundary condition}.
\InkRed{Other boundary conditions can be either converted to the homogeneous Dirichlet boundary condition \cite{Ern04} or treated without major modifications.}
The corresponding variational problem is to find  a minimizer $u \in U_0$ such that
\begin{align}
   I[u]=\min_{v\in U_0}I[v].
   \label{min-1}
\end{align}
\InkRed{We assume that $F$ satisfies suitable conditions so that the variational problem \eqref{min-1} has a unique solution.
We also assume that a linearization of the first variation of the functional (i.e. $B[u;\cdot,\cdot]$, see \eqref{eq:Bu} below) is convex around the solution of the minimization problem.
The latter is needed for the error problem \eqref{eq:error-problem} defined in Section~\ref{ssec:hbee}
to have a unique solution and for its Gau{\ss}-Seidel solution to converge.
However, we do not assume that $F$ is quadratic
although some theoretical analysis in this paper is valid only for convex quadratic functionals.}

A necessary condition for $u$ to be a minimizer is that the first variation of the functional vanishes.
This leads to the Galerkin formulation 
\begin{align*}
   \delta I [u,v] 
   \equiv \int_{\Omega} \left( F_{u}(\boldx,u,\nabla u) \; v 
         + F_{\nabla u}(\boldx,u,\nabla u)\cdot \nabla v\right) d\boldx 
   = 0,
   \qquad \forall v \in U_{0}
   % \label{galerkin-1}
\end{align*}
where $F_{u}$ and $F_{\nabla u}$ are the partial derivatives of $F$ with respect to $u$ and $\nabla u$, respectively.
The linearization of the first variation $\delta I[\cdot,\cdot]$ around the minimizer $u$ with respect to the first argument
is given by
\begin{align}
  B[u;w,v] = & \int_{\Omega} \left (\frac{}{}
       v \frac{\partial^2 F}{\partial u^2}(\boldx,u,\nabla u) w
      + v \frac{\partial^2 F}{\partial u \, \partial \nabla u} (\boldx,u,\nabla u) \cdot \nabla w 
      +  w \frac{\partial^2 F}{\partial u \, \partial \nabla u} (\boldx,u,\nabla u) \cdot \nabla v
      \right .
     \notag \\
      & \left.  
      + \nabla w\cdot \frac{\partial^2 F}{\partial \nabla u \partial \nabla u}(\boldx,u,\nabla u) \cdot \nabla v \frac{}{}
      \right ) d \boldx ,\qquad \forall w, v \in U_{0} .
   \label{eq:Bu}
\end{align}

%--- Adaptive linear finite element approximation ------------------------
Given a triangulation $\cT_h$ for $\Omega$, let $U_0^h \subset U_0$ be the associated linear finite element space.
Then, a linear finite element approximation of $u$ is defined as $u_h \in U_0^h$ such that
\begin{align*}
   I[u_h]=\min_{v_h\in U_g^h}I[v_h].
   % \label{min-1h}
\end{align*}
The solution $u_h$ can also be found by solving the corresponding Galerkin formulation
\begin{align}
   \delta I [u_h, v_h] = \int_\Omega \left(F_u(\boldx,u_h,\nabla u_h) \; v_h 
      + F_{\nabla u}(\boldx,u_h,\nabla u_h) \cdot \nabla v_h \right) d\boldx
      = 0,
      \qquad \forall v_h \in U_0^h .
   \label{fem-2}
\end{align}
%--- mesh adaptation ---
For the adaptive finite element solution, the mesh $\cT_h$ is generated according to the behaviour of the error of the approximation $u_h$.
We follow  the so-called $M$-uniform mesh approach \cite{Huang05a} with which an adaptive mesh is generated as a uniform mesh in the metric specified by a symmetric and strictly positive definite tensor $M = M(\boldx)$. 
Such a mesh is called an \emph{$M$-uniform mesh}.
A scalar metric tensor will lead to an isotropic mesh while a full metric tensor will generally result in an anisotropic mesh.
In this sense, the mesh generation procedure is the same for both isotropic and anisotropic mesh generation.
The key to the approach is to choose a proper metric tensor $M$.

Once a metric tensor $M$ has been chosen, a sequence of mesh and corresponding finite element approximation are generated in an iterative fashion.
We start with an initial mesh $\cT_h^{(0)}$.
For every mesh $\cT_h^{(i)}$ we solve the variational problem with $U_{0,(i)}^h$ for finite element approximation $u_h^{(i)}$.
The new mesh $\cT_h^{(i+1)}$ is generated as an almost $M$-uniform mesh in the metric specified by $M_h^{(i)}$ which is computed based on $u_h^{(i)}$ and $\cT_h^{(i)}$.
This yields the sequence
\[
   (\cT_h^{(0)}, U^h_{0,(0)}) \rightarrow u_h^{(0)} \rightarrow M_h^{(0)}
   \rightarrow 
   (\cT_h^{(1)}, U^h_{0,(1)}) \rightarrow u_h^{(1)} \rightarrow M_h^{(1)}
   \rightarrow \cdots 
\]
The process is repeated until a good adaptation is achieved.

In the following section we derive the metric tensor $M$ for use in anisotropic mesh adaptation for the variational problem \eqref{min-1}.  
First, we derive a bound on the variation $\delta I [u_h, v]$ and then define $M$ such that the obtained bound is minimized on corresponding $M$-uniform meshes.

%**************************************************************************
\section{Monitor functions based on hierarchical basis error estimation}
%**************************************************************************
\label{sec:HBEE}

\subsection{Error bounds for variational problems}
%--------------------------------------------------------------------------

%--- motivation ---
Let $e_h = u - u_h$ be the error of the finite element solution $u_h$.
It is known  \cite{Evans98} that for a uniformly convex quadratic functional in the form \eqref{functional-1}, there exists a positive constant $\beta$ such that
\begin{align*}
   \| \nabla e_h \|_{L^2(\Omega)}^2 
      \le \beta \left| \delta I[e_h,e_h] \right| 
      = \beta \left| \delta I[u_h,e_h] \right|.
   % \label{bound-i4}
\end{align*}
In this case, the quantities
\begin{align}
   \left|\delta I[u_h,e_h]\right| \quad \mbox{and}\quad
   \left(\frac{|\delta I[u_h,e_h]|}{\| \nabla e_h \|_{L^2(\Omega)}}\right)^2
   \label{variation-1}
\end{align}
are equivalent to $\| \nabla e_h \|_{L^2(\Omega)}^2$, and minimizing their bounds is equivalent to minimizing error bounds.
Consequently, it is reasonable to define $M$ based on bounds for these quantities even for other functionals.

%--- derivation ---
We begin with deriving a bound for $|\delta I[u_h, e_h]|$.
For simplicity we denote
\[   F_u(\boldx) 
      = F_u(\boldx,u_h(\boldx),\nabla u_h(\boldx)),\quad F_{\nabla u}(\boldx)
      = F_{\nabla u}(\boldx, u_h(\boldx),\nabla u_h(\boldx)).
\]
Then \InkRed{for any $v \in U_0$}
\begin{align}
   \delta I [u_h, v] 
   & = \int_\Omega \left [F_u(\boldx,u_h,\nabla u_h) \; v 
      + F_{\nabla u}(\boldx,u_h,\nabla u_h)
         \cdot \nabla v \right ] d \boldx \notag \\
   & = \sum_{K \in \mathcal{T}_h} \int_{K} \left [ 
      F_u(\boldx ,u_h,\nabla u_h) \; v 
         + F_{\nabla u}(\boldx,u_h,\nabla u_h)
         \cdot \nabla v \right ] d \boldx \notag \\
   & = \sum_{K \in \mathcal{T}_h} \int_{K} \left [ 
      F_u(\boldx) \; v 
         + F_{\nabla u}(\boldx) \cdot \nabla v \right ] d \boldx .
   \label{disc-1}
\end{align}
From the divergence theorem, \eqref{disc-1} can be rewritten as
\begin{align}
   \delta I [u_h, v] 
      = \sum_{K \in \cT_h} \int_{K} \left[ F_u(\boldx) 
         - \nabla \cdot F_{\nabla u}(\boldx) \right] \; v \, d\boldx 
      + \sum_{K \in \cT_h} \int_{\partial{K}} v \;
         F_{\nabla u}(\boldx) \cdot \boldsymbol{n} \, ds, 
   \label{disc-2}
\end{align}
where $\boldsymbol{n}$ denotes the outward unit normal to the face $\partial K$.

Now, let $\partial \cT_h$ be the collection of all faces of mesh $\cT_h$ and $K$ and $K'$ be the two elements sharing a common face $\gamma$.
We define residual $r_h$ and edge jump $R_h$ as
\begin{align*}
   r_{h}(\boldx) & 
      = F_{u}(\boldx)-\nabla \cdot F_{\nabla u}( \boldx) 
         \qquad \forall \boldx \in K \quad \forall K \in \cT_h,
   % \label{res-int} 
   \\
   R_{h}(\boldx) &= 
   \begin{cases}
   (F_{\nabla u}(\boldx)\cdot \boldsymbol{n}_\gamma)|_{K}
      +(F_{\nabla u}(\boldx)\cdot\boldsymbol{n}_\gamma)|_{K'}
      & \quad \boldx \in \gamma 
         \quad \forall \gamma \in \partial \cT_h\backslash \partial \Omega,\\ 
   0 & \quad \boldx \in \gamma \quad \forall \gamma \in \partial \Omega.
   \end{cases}
   % \label{res-bnd}
\end{align*}
Then \eqref{disc-2} becomes 
\begin{align*}
   \delta I[u_{h},v] 
      = \sum_{K \in \cT_h} \int_{K} r_{h}(\boldx) \; v \, d\boldx 
         + \sum_{\gamma \in \partial \cT_h}
            \int_\gamma R_h(\boldx) \; v\, ds.
\end{align*}
Taking $v = e_h = u - u_h$ in the above equation and using Schwarz's inequality, we obtain
\begin{align}
  \abs{\delta I[u_h, e_h]}
   &\le \sum_{K \in \cT_h} \int_K \abs{r_h(\boldx) \; e_h} d\boldx
      + \sum_{\gamma \in \partial \cT_h} 
         \int_\gamma \abs{R_h(\boldx) \; e_h} ds \notag \\
   & \le \sum_{K \in \cT_h} \norm{r_h}_{L^2(K)} \norm{e_h}_{L^2(K)}
      + \sum_{\gamma \in \partial \cT_h} \norm{R_h}_{L^2(\gamma)} \,
         \norm{e_h}_{L^2(\gamma)} 
   % \label{bound-1a} 
   \notag \\
   &= \sum_{K \in \cT_h} \left[\norm{r_h}_{L^2(K)} \norm{e_h}_{L^2(K)}
      + \frac{1}{2} \sum_{\gamma \in \partial K} \norm{R_h}_{L^2(\gamma)}
         \norm{e_h}_{L^2(\gamma)} \right].
   \label{bound-1b}
\end{align}
For further derivation, we need to estimate $\norm{e_h}_{L^2(K)}$ and $\norm{e_h}_{L^2(\gamma)}$. 
In \cite{HuaLi10}, the authors employ anisotropic interpolation error bounds on elements and element faces.
This results in a semi-a~posteriori metric tensor which involves residual and edge jumps, both dependent on the computed solution, and the Hessian of the exact solution.
In the following, we use an a~posteriori hierarchical basis error estimate.

\subsection{An a~posteriori error estimator based on hierarchical bases}
%--------------------------------------------------------------------------
\label{ssec:hbee}
The computation of the error estimator is based on a general framework, details on which can be found among others in the work of Bank and Smith \cite{BanSmi93} or Deuflhard et al. \cite{DeLeYs89}.
The approach is briefly explained as follows.

Recall that $u_h \in U_0^h$ is a linear finite element solution of the Galerkin formulation \eqref{fem-2} and the error is $e_h = u - u_h$.
\InkRed{Consider the extension of $U_0^h$ to a subspace of piecewise quadratic functions with the linear span of the edge bubble functions $W^h$, i.e.
\[
   U_0^h \rightarrow U_0^h \oplus W^h.
\]
}
Recall also that 
\begin{align*}
   \delta I[u_h + e_h, v] = 0 \qquad \forall v \in U_0.
\end{align*}
%Denote by $B[u_h;\cdot,\cdot]$ a bilinear form resulting from a linearization
%of $\delta I[\cdot,\cdot]$ about $u_h$ with respect to the first argument.
The error estimate $z_h$ is then defined as the solution of the approximate \emph{linear} error problem 
\begin{align}
   \begin{cases}
      \text{Find } z_h\in W^h \text{ such that} \\
      \delta I[u_h,w_h] + B[u_h;z_h,w_h] = 0 \quad \forall w_h \in W^h .
   \end{cases}
   \label{eq:error-problem}
\end{align}
Recall that $B[u_h;\cdot,\cdot]$ is the bilinear form resulting from the linearization
of $\delta I[\cdot,\cdot]$ about $u_h$ with respect to the first argument.
The estimate $z_h$ can be viewed as a projection of the error onto the subspace $W^h$. 
\InkRed{By construction, $\Pi_h z_h = 0$ for the linear finite element interpolation operator $\Pi_h$.
Thus, $z_h - \Pi_h z_h = z_h$.}

This definition of the error estimate is global and its solution can be costly.
To avoid the expensive exact solution in numerical computation, we employ only a few sweeps of the symmetric
Gau{\ss}-Seidel iteration for the resulting linear system, which proves to be sufficient for the purpose of mesh adaptation (see \cite{HuKaLa10} or numerical examples in Section~\ref{sec:Numerical-Results}).

\InkRed{It is noted that under the convex assumption of $B[u;\cdot,\cdot]$, the error problem (\ref{eq:error-problem})
is solvable and has a unique solution when $u_h$ is sufficiently close to $u$. Moreover, Gau{\ss}-Seidel iteration
is convergent when applied to the solution of (\ref{eq:error-problem}).}

\subsection{Optimal metric}
%--------------------------------------------------------------------------
\InkRed{We now use the error estimate $z_h$ to develop the metric tensor.
Instead of directly using the bound \eqref{bound-1b}, we use
\begin{align}
   E[u_{h},z_h] = \sum_{K \in \cT_h} \left[\norm{r_h}_{L^2(K)} \norm{z_h}_{L^2(K)}
      + \frac{1}{2} \sum_{\gamma \in \partial K} \norm{R_h}_{L^2(\gamma)}
      \norm{z_h}_{L^2(\gamma)} \right] ,
   \label{eq:e-u-z}
\end{align}
which is obtained by replacing $e_h$ in \eqref{bound-1b} with $z_h$.
Theoretically, it is unclear if the bound \eqref{bound-1b} is bounded by $E[u_h,z_h]$
although our numerical results show that the metric tensor based on $E[u_h,z_h]$ leads to
proper mesh adaptation for all examples we have tested.
On the other hand, this is somewhat related to the saturation assumption, which basically states that the quadratic finite element approximation provides a truly better approximation than the linear one.
The saturation assumption is known to be valid for some cases
under suitable conditions for isotropic meshes \cite{DoeNoc02}.
It is still unknown if such results hold for anisotropic meshes.
}
% We now assume that $z_h$ provides a reliable local estimate on $e_h$, i.e., there exist constants $C_1 > 0$ and $C_2 > 0$ such that
% \[
%    \norm{e_h}_{L^2(K)} \le C_1 \norm{z_h}_{L^2(K)}
%       \quad \text{and} \quad
%    \norm{e_h}_{L^2(\gamma)} \le C_2 \norm{z_h}_{L^2(\gamma)} .
% \]
% Then we can replace $e_h$ with $z_h$ in \eqref{bound-1b} and develop the bound on $\delta I[u_{h},e_h]$ in terms of $z_h$.

Recalling that $\Pi_h z_h = 0$ and using element-wise interpolation error estimates in \cite{Huang05a}, we have
\begin{align}
   \norm{z_h}_{L^{2}(K)}
   & = \norm{z_h - \Pi_h z_h}_{L^{2}(K)} \notag \\
   & \leq C\left[ \int_K 
     \left( \tr \left( (F_K')^T \abs{H(z_h)} F_K'\right) \right)^2 d\boldx
         \right]^{\frac{1}{2}} \notag \\
   & = C \abs{K}^\frac{1}{2} \tr\left( (F_K')^T \abs{H_K(z_h)} F_K'\right),
   \label{bound-2} \\
   \left[ \sum \limits_{\gamma \in \partial K}\frac{1}{\abs{\gamma}}
      \norm{z_h}_{L^2(\gamma)}^2 \right]^{\frac{1}{2}}
   & = \left[ \sum \limits_{\gamma \in \partial K}\frac{1}{\abs{\gamma}}
   \norm{z_h - \Pi_h z_h}_{L^2(\gamma)}^2 \right]^{\frac{1}{2}} \notag \\
      & \leq C \left[ \frac{1}{\abs{K}} \int_K 
      \left( \tr\left( (F_K')^T \abs{H(z_h)} F_K' \right)\right)^2 d\boldx
      \right]^{\frac{1}{2}} \notag \\
      & = C  \tr\left( (F_K')^T \abs{H_K(z_h)} F_K' \right) ,
   \notag
\end{align}
where $\tr(\cdot )$ denote the trace of a matrix, $H_K(z_h)$ is the Hessian of $z_h$ on $K$, $\abs{H_K(z_h)}= \sqrt{H_K^2(z_h)}$, $\abs{K}$ is the volume of $K$, $F_K$ is a mapping from the reference element $\hat K$ to element $K$, and $C$ is a constant independent of $\cT_h$ and $z_h$.
Thus,
\begin{align}
    \sum_{\gamma\in\partial K} \norm{R_h}_{L^2(\gamma)} 
      \norm{z_h}_{L^2(\gamma)} 
   & = \sum_{\gamma\in\partial K} \abs{\gamma}^{\frac{1}{2}}
      \norm{R_h}_{L^2(\gamma)} \abs{\gamma}^{-\frac{1}{2}} 
         \norm{z_h}_{L^2(\gamma)} \notag \\
   & \leq \left(\sum_{\gamma\in\partial K} \abs{\gamma}
         \norm{R_h}_{L^2(\gamma)}^2\right)^{\frac{1}{2}}
      \left(\sum_{\gamma\in\partial K}\frac{1}{\abs{\gamma}}
         \norm{z_h}_{L^2(\gamma)}^2\right)^{\frac{1}{2}} \notag \\
   & \leq \left(\sum_{\gamma\in\partial K} \abs{\gamma}^{\frac{1}{2}}
      \norm{R_h}_{L^2(\gamma)}\right)
         \left(\sum_{\gamma\in\partial K}\frac{1}{\abs{\gamma}} 
            \norm{z_h}_{L^2(\gamma)}^2 \right)^{\frac{1}{2}}
      \notag \\
   & \leq C \left( \sum_{\gamma\in\partial K}\abs{\gamma}^{\frac{1}{2}}
         \norm{R_h}_{L^2(\gamma)}\right) 
      \; \tr\left( (F_K')^T \abs{H_K(z_h)} F_K' \right).
   \label{bound-3c}
\end{align}
Substituting \eqref{bound-2} and \eqref{bound-3c} into \eqref{eq:e-u-z} leads to
\begin{align}
   E[u_h,z_h] 
   & \leq C \sum_{K\in\cT_h} \left[\abs{K}^{\frac{1}{2}}\norm{r_h}_{L^2(K)}
         + \sum_{\gamma\in\partial K}\abs{\gamma}^{\frac{1}{2}}
            \norm{R_h}_{L^2(\gamma)} \right]
         \; \tr\left( (F_K')^T \abs{H_K(z_h)} F_K' \right) \notag \\
      & =  C \sum_{K\in\cT_h}
         \left[\frac{\norm{r_h}_{L^2(K)}}{\abs{K}^{\frac{1}{2}}}
            + \frac{1}{\abs{K}} 
               \sum_{\gamma\in\partial K}\abs{\gamma}
               \frac{\norm{R_h}_{L^2(\gamma)}}{\abs{\gamma}^{\frac{1}{2}}}
         \right]
       \; \abs{K} \tr\left( (F_K')^T \abs{H_K(z_h)} F_K' \right) \notag \\
      & =  C \sum_{K\in\cT_h}
         \left[\snorm{r_h}_{L^2(K)}
         + \frac{1}{\abs{K}} \sum_{\gamma\in\partial K}\abs{\gamma}
            \snorm{R_h}_{L^2(\gamma)} \right]
       \; \abs{K} \tr\left( (F_K')^T \abs{H_K(z_h)} F_K' \right),
   \label{bound-5}
\end{align}
where $\snorm{\cdot}$ denotes the scaled norm
\[
   \snorm{v}_{L^2(\omega)} 
      = \frac{1}{\abs{\omega}^{1/2}} \norm{v}_{L^2(\omega)}
      = \left[
            \frac{1}{\abs{\omega}}\int_\omega v^2 d\boldx
        \right]^{\frac{1}{2}}.
\]

We now use bound (\ref{bound-5}) to define the metric tensor $M$. 
To ensure that $M$ is strictly positive definite, we first regularize the bound with a positive constant $\alpha _{h}$ (to be determined), i.e., 
\begin{align}
   E[u_h,z_h] 
      & \le C \sum_{K\in \cT_h} \left[ \alpha_h + \snorm{r_h}_{L^2(K)}
         + \frac{1}{\abs{K}} \sum_{\gamma\in\partial K} \abs{\gamma}
            \snorm{R_h}_{L^2(\gamma)} \right] \notag \\
      & \qquad \qquad \times \abs{K} \tr\left( (F_K')^T \left(\alpha_h I 
            + \abs{H_K(z_h)}\right) F_K' \right) \notag \\
      &  = C \alpha_{h}^{2} \sum_{K\in \cT_h} \left[ 
         1 + \frac{1}{\alpha_h} \snorm{r_h}_{L^2(K)}
         + \frac{1}{\alpha_h \abs{K}}\sum_{\gamma\in\partial K}\abs{\gamma}
               \snorm{R_h}_{L^2(\gamma)} \right] \notag \\
      & \qquad \qquad 
         \times \abs{K} \tr\left( (F_K')^T H_{K,\alpha}(z_h) F_K' \right),
   \label{bound-6}
\end{align}
where
\[
   H_{K,\alpha}(z_h) = I + \frac{1}{\alpha_h} \abs{H_K(z_h)}.
\]

The optimal metric tensor is obtained by minimizing bound (\ref{bound-6}) for $M$-uniform meshes.
It is known \cite{Huang06} that an $M$-uniform mesh satisfies the \emph{alignment condition}
\begin{align}
   \frac{1}{d} \tr \left((F_K')^T M_K F_K'\right)
      & = \det \left((F_K')^T M_K F_K'\right)^{\frac{1}{d}}
   \label{align-1}
\end{align}
and the \emph{equidistribution condition}
\begin{align}
   \rho_K\abs {K} 
      & = \frac{\sigma_h}{N} ,
   \label{eqdis-1}
\end{align}
where $\sigma_h = \sum_{K\in \cT_h} \rho_K \abs{K}$, $M_K$ is an average of $M$ over $K$, $\rho_K = \sqrt{\det (M_K)}$, and $N$ the number of elements of $\cT_h$.

We now pay our attention to the $\mbox{tr}(\cdot)$ factor in (\ref{bound-6}).
Notice that in general,
\[
   \frac{1}{d} \tr \left((F_K')^T H_{K,\alpha}(z_h) F_K'\right)
   \geq
   \det \left((F_K')^T H_{K,\alpha}(z_h) F_K'\right)^{\frac{1}{d}}.
\]
From (\ref{align-1}) we can see that the equality in the above inequalities holds if we choose $M = M_K$ in the form 
\begin{align}
   M_K = \theta_K H_{K,\alpha}(z_h) \qquad \forall K \in \cT_h
   \label{metric-an}
\end{align}
for some scalar function $\theta =\theta_K$.
Indeed, with this choice of $M_K$ the alignment condition \eqref{align-1} reads as
\begin{align}
   \frac{1}{d} \tr \left((F_K')^T H_{K,\alpha} F_K'\right) 
      = \det \left((F_K')^T H_{K,\alpha } F_K'\right)^{\frac{1}{d}}
      = \abs{K}^{\frac{2}{d}} 
         \det \left( H_{K,\alpha}(z_h) \right)^{\frac{1}{d}} ,
   \label{eqdis-1b}
\end{align}
where we have used $\abs{\det (F_K')} = \abs{K}$ and assumed $|\hat K| = 1$.
Substituting \eqref{eqdis-1b} into \eqref{bound-6} yields
\begin{align}
   E[u_h,z_h]
      & \le C \alpha_{h}^{2} \sum_{K\in \cT_h} \left[ 
         1 + \frac{1}{\alpha_h} \snorm{r_h}_{L^2(K)}
         + \frac{1}{\alpha_h\abs{K}} \sum_{\gamma\in\partial K}\abs{\gamma}
               \snorm{R_h}_{L^2(\gamma)} \right] \notag \\
      & \qquad \qquad \times \abs{K}^{\frac{d+2}{d}} 
         \det \left( H_{K,\alpha}(z_h) \right)^{\frac{1}{d}}. 
   \label{bound-an}
\end{align}
Next, we use the equidistribution condition \eqref{eqdis-1} to determine $\theta =\theta_K$ in \eqref{metric-an}.
From H{\"o}lder's inequality, we have 
\begin{align*}
   \left[ \frac{1}{N}
      \sum_K \left(\abs{K}\rho_K\right)^{\frac{d+2}{d}}\right]^{\frac{d}{d+2}}
   \geq \frac{1}{N}\sum_K \abs{K}\rho_K = \frac{\sigma_h}{N}
\end{align*}
or
\begin{equation}
   \sum_{K}\left( |K|\rho _{K}\right) ^{\frac{d+2}{d}}\geq
   \left( {\sigma _{h}}\right) ^{\frac{d+2}{d}}N^{-\frac{2}{d}} ,
   \label{eqdis-2}
\end{equation}
with the lower bound being attained for a mesh satisfying (\ref{eqdis-1}).
Comparing the left-hand side of \eqref{eqdis-2} with the right-hand side of \eqref{bound-an} suggests that $\rho =\rho_K$ be defined as
\begin{align}
   \rho_K 
      & = \left[1 + \frac{1}{\alpha_h}\snorm{r_h}_{L^2(K)}
      + \frac{1}{\alpha_h\abs{K}} \sum_{\gamma\in\partial K} \abs{\gamma}
            \snorm{R_h}_{L^2(\gamma)}\right]^{\frac{d}{d+2}}
            \det \left(H_{K,\alpha}(z_h)\right)^{\frac{1}{d+2}}  \notag \\
      & = \left[1 + \frac{1}{\alpha_h}\snorm{r_h}_{L^2(K)}
      + \frac{1}{\alpha_h\abs{K}} \sum_{\gamma\in\partial K} \abs{\gamma}
            \snorm{R_h}_{L^2(\gamma)}\right]^{\frac{d}{d+2}}
       \det 
       \left(I + \frac{1}{\alpha_h}\abs{H_K(z_h)}\right)^{\frac{1}{d+2}} .
   \label{rho-an1}
\end{align}
From relations $\rho_K =\sqrt{\det (M_K)}$ and $M_K = \theta_K H_{K,\alpha}(z_h)$ we can obtain $\theta_K$.
The metric tensor $M_K$ is then given by
\begin{align}
   M_K = \rho_K^{\frac{2}{d}}
      \det \left(I + \frac{1}{\alpha_h}\abs{H_K(z_h)}\right)^{-\frac{1}{d}}
      \left[I + \frac{1}{\alpha_h}\abs{H_K(z_h)}\right] .
   \label{metric-an1}
\end{align} 

With this choice of $\rho_K$ (and $M_K$), the right-hand side of (\ref{bound-an}) attains its lower bound for a mesh satisfying the equidistribution condition (\ref{eqdis-1}).
Then, the variation of the functional \eqref{functional-1} has an upper bound as
\begin{align*}
   E[u_h,z_h]
      \le C \alpha_h^2 (\sigma_h)^{\frac{d+2}{d}} N^{-\frac{2}{d}}.
\end{align*}

To complete the definition of the metric tensor, we need to choose the regularity parameter $\alpha_h$.
Following \cite{Huang05a}, we choose it such that
\[
   \sigma_h 
      \equiv \sum_{K\in\cT_h}\rho_K \abs{K} 
      \boldsymbol{=} 2\abs{\Omega}
\]
or
\begin{align}
   &\sum_{K\in\cT_h}\abs{K} 
      \left[1 + \frac{1}{\alpha_h}\snorm{r_h}_{L^2(K)}
         + \frac{1}{\alpha_h \abs{K}}\sum_{\gamma \in \partial K}
            \abs{\gamma}\snorm{R_h}_{L^2(\gamma)}\right]^{\frac{d}{d+2}}
      \det \left(I + \frac{1}{\alpha_h}\abs{H_K(z_h)}\right)^{\frac{1}{d+2}}
         \notag \\
      & \qquad = 2 \abs{\Omega}. 
   \label{alpha-an}
\end{align}
With this choice, roughly fifty percents of the mesh elements will be concentrated in the regions of large $\rho$ \cite{Huang05a}. 
It is easy to show that \eqref{alpha-an} has a unique solution since its left-hand side is monotonically decreasing with $\alpha_h$ increasing and tends to $\abs{\Omega}$ as $\alpha_h \to \infty$ and to $+\infty$ as $\alpha_h \to 0$.
Moreover, it can be solved using a simple iteration scheme such as the bisection method.

Once $\alpha_h$ is computed, the adaptation function $\rho$ and the metric tensor $M$ can be determined by \eqref{rho-an1} and \eqref{metric-an1}, respectively.
A new mesh can then be generated based on the metric tensor.

This definition of the metric tensor involves the residual $r_h$, the edge jump $R_h$, and the HBEE $z_h$.
All these quantities are based on the computed solution; in this sense, (\ref{metric-an1}) is a~posteriori. 

%**************************************************************************
\section{Numerical results}
%**************************************************************************
\label{sec:Numerical-Results}
In this section, we present some numerical results for a selection of two-dimensional problems with an anisotropic behaviour and compare the following metric tensors:
\begin{itemize}
  \item isotropic a~posteriori metric tensor \cite{HuaLi10} based on the variational formulation, residual, and edge jumps,
  \item anisotropic semi-a~posteriori metric tensor \cite{HuaLi10} based on the variational formulation, residual, edge jumps, and Hessian recovery,
  \item anisotropic a~posteriori metric tensor based on the variational formulation, residual, edge jumps, and hierarchical basis error estimator (HBEE) developed in the previous section.
\end{itemize}

The HBEE is computed as a global error problem.
A few symmetric Gauß--Seidel iterations are employed to obtain an approximation for $z_h$ until the relative difference of the old and the new approximations is under a given tolerance \emph{GS-RTOL} $= 0.01$.
More details on several implementation issues such as mesh quality measure and computation of the error estimator can be found in \cite[Section~4]{HuKaLa10}.

\subsection{A quadratic functional}
\label{ex:tanh}
%--------------------------------------------------------------------------
% tanh example: convergence plots
\begin{figure}[t] \centering
   \subfloat[Example~\ref{ex:tanh}: the exact solution.]{
      \includegraphics[width=0.38\textwidth,clip]{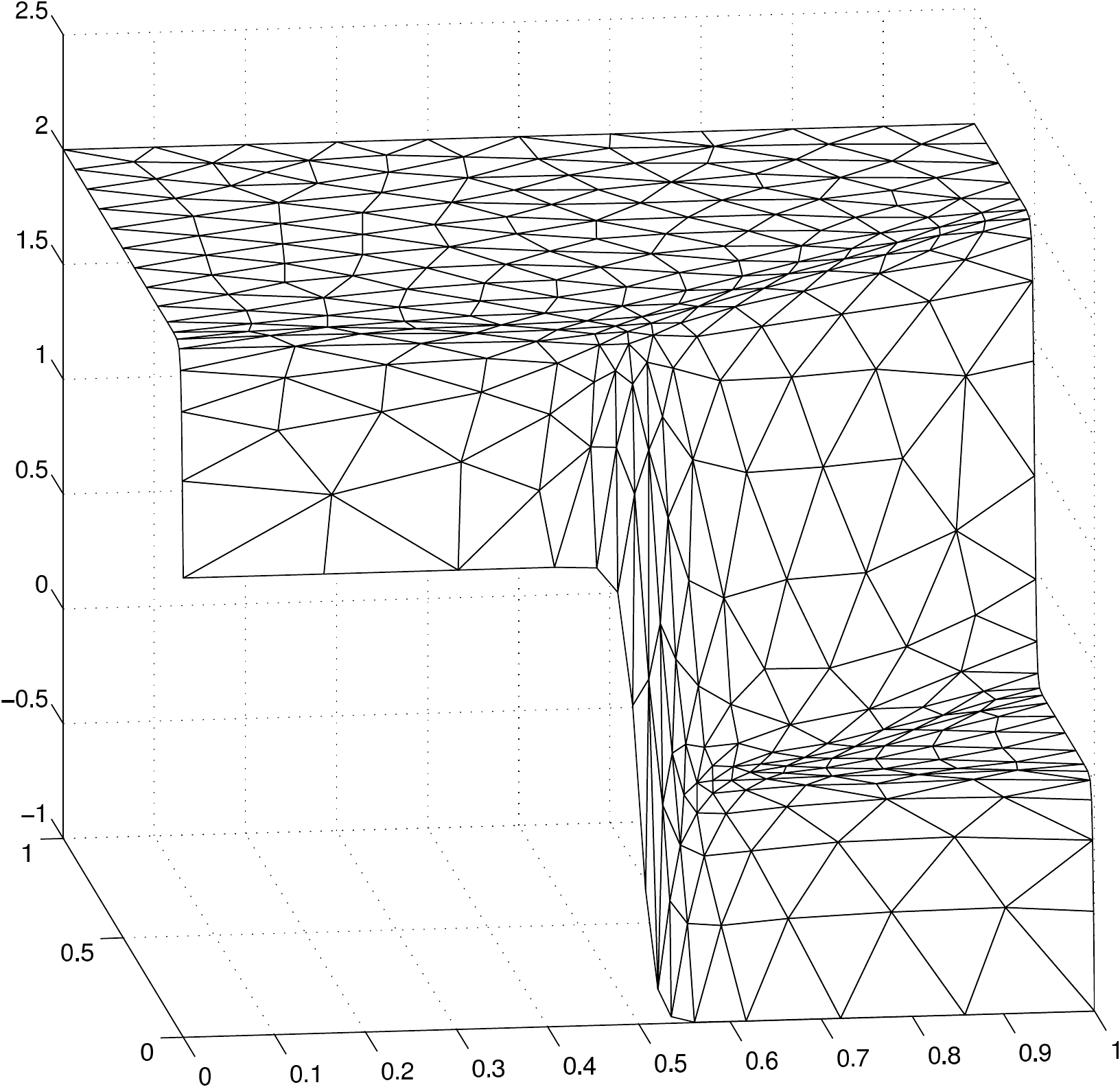}
      \label{fig:tanh:plot}
   } \qquad
   \subfloat[Example~\ref{ex:tanh}: error comparison.]{
      \includegraphics[width=0.50\textwidth,clip]{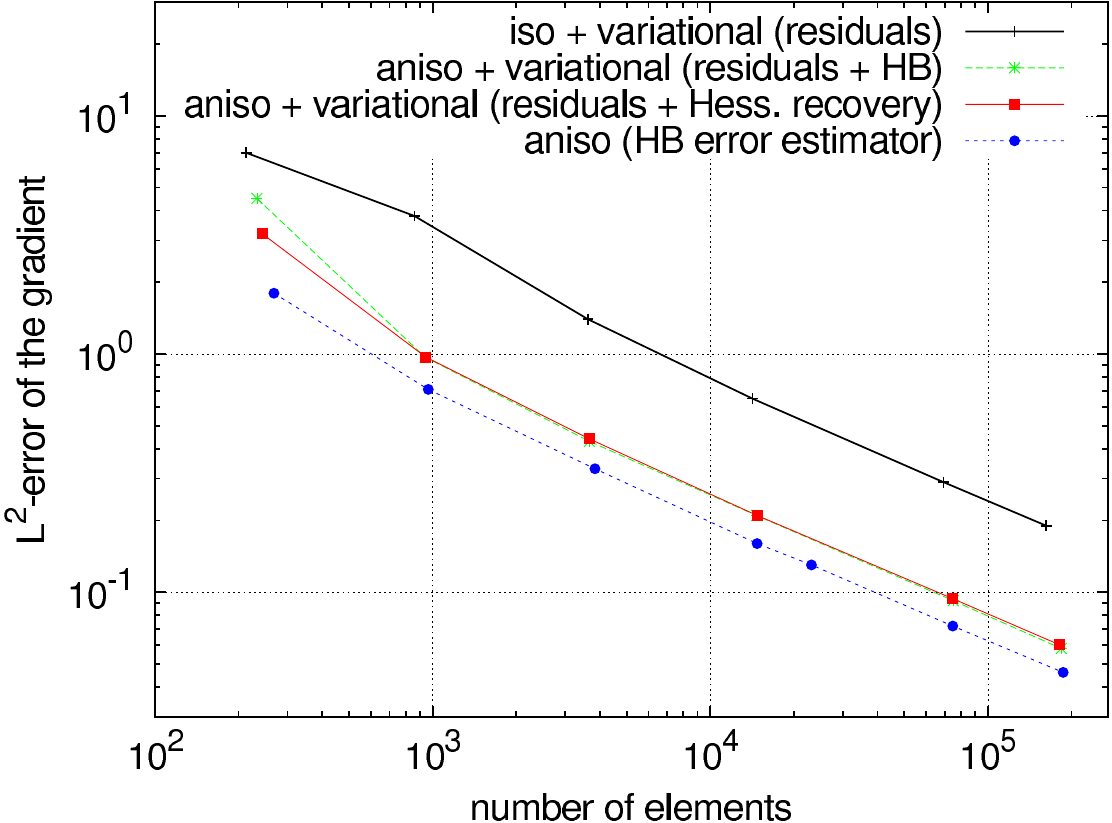}
     \label{fig:tanh:convergence}
   }
   \caption{Example~\ref{ex:tanh}: 
      \protect \subref{fig:tanh:plot} surface plot of the exact solution,
      \protect \subref{fig:tanh:convergence} the finite element error 
      $\norm{\nabla e_h}_{L^2(\Omega)}$ against the number of elements.
   }
\end{figure}

As a first example we consider a quadratic functional 
\begin{align}
   I[u] = \int_\Omega \left[ \frac{1}{2} \abs{\nabla u}^2 - u f \right] d\boldx
\label{eq:bvp-vp}
\end{align}
with $\Omega = (0,1)\times(0,1)$.
The function $f$ and the Dirichlet boundary conditions are chosen such that the exact solution of the minimization problem \eqref{min-1} is given by
\begin{align*}
   u_{exact}(x,y) = \tanh(60x) - \tanh\left(60(x - y) - 30\right).
\end{align*}
A surface plot of $u_{exact}$ is given in Fig.~\ref{fig:tanh:plot}.
The solution exhibits a strong anisotropic behaviour and describes the interaction between a boundary layer along the $x$-axis and a shock wave along the line $y = x-0.5$.

Note that this variational problem is equivalent to the boundary value problem
\begin{align*}
   \begin{cases}
     -\Delta  u = f  & \text{in } \Omega, \\
              u = u_{exact}  & \text{on } \partial \Omega.
   \end{cases}
\end{align*}
Thus, we  also compare metric tensors based on the variational formulation described above with the one based on the hierarchical basis error estimator and PDE developed in \cite{HuKaLa10}. 

\begin{figure}[t] \centering
      \includegraphics[width=0.30\textwidth,clip]{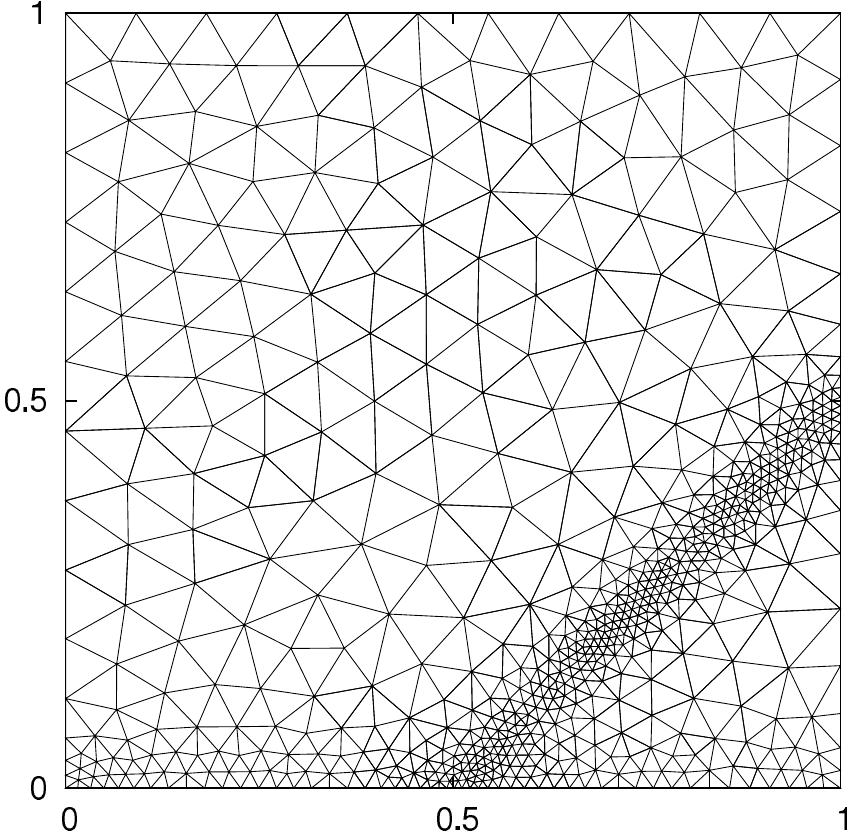}
      \qquad
      \includegraphics[width=0.31\textwidth,clip]{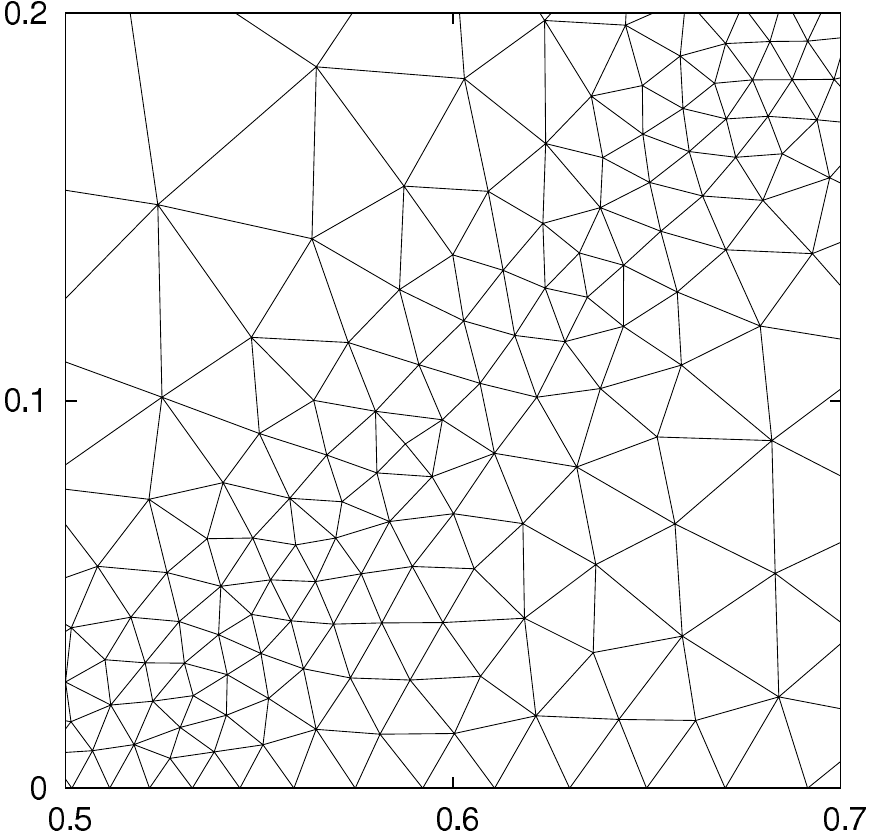}
      \caption{Example~~\ref{ex:tanh} (isotropic): 
         adaptive mesh and a close-up view
         at (0.5,0) obtained with the isotropic metric tensor based 
         on variational formulation with residual and edge jumps; 
         664 vertices and 1226 triangles, 
         $\norm{\nabla e_h}_{L^2} = 2.8$, maximum aspect ratio 3.}
      \label{fig:tanh:iso}
\end{figure}

\begin{figure}[p] \centering
    \subfloat[Variational formulation with residual, edge jumps, and Hessian recovery:
         663 vertices and 1235 triangles, 
         $\norm{\nabla e_h}_{L^2} = 6.2 \times 10^{-1}$,
         maximum aspect ratio 34.]{
      \includegraphics[width=0.30\textwidth,clip]{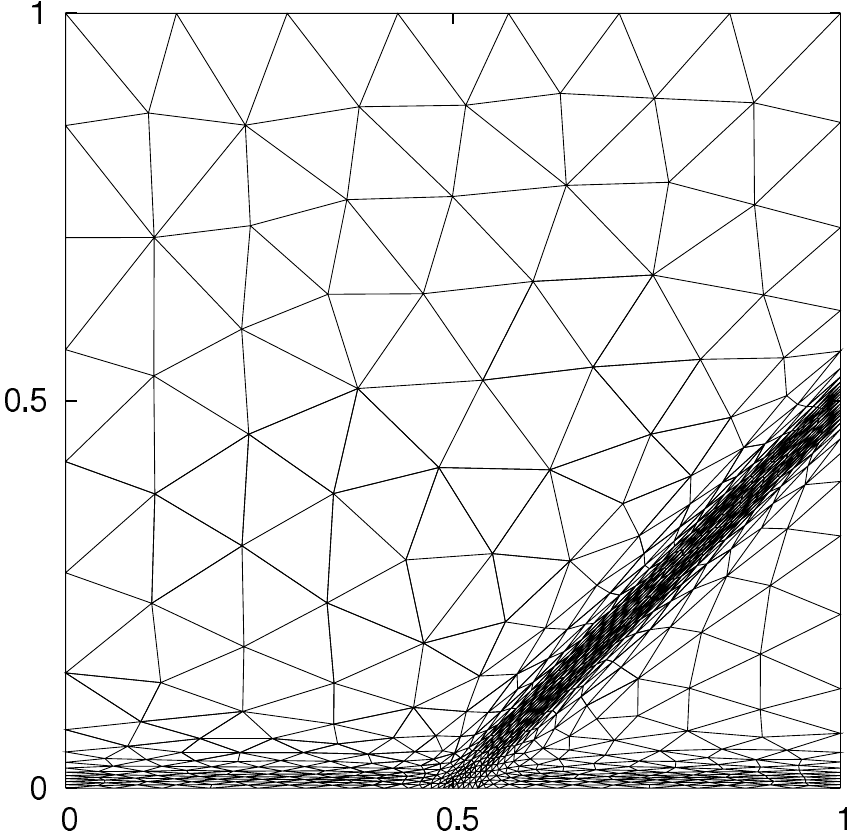}
      \qquad
      \includegraphics[width=0.31\textwidth,clip]{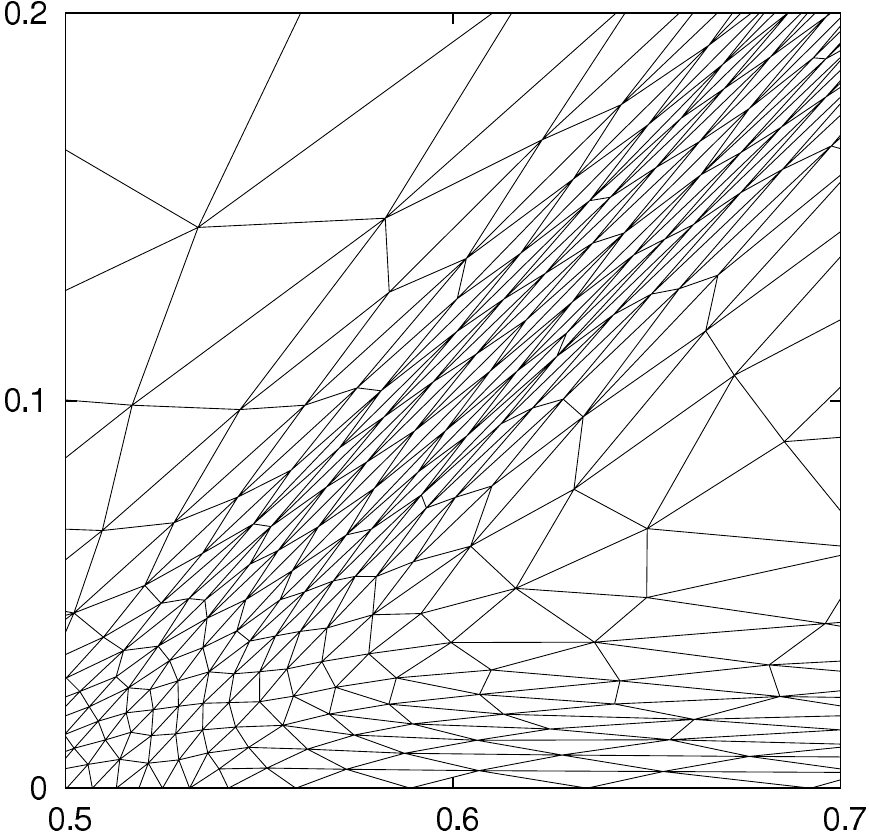}
   }\\
   \subfloat[Variational formulation with residual, edge jumps, and HBEE:
         672 vertices and 1249 triangles,
         $\norm{\nabla e_h}_{L^2} = 6.0 \times 10^{-1}$,
         maximum aspect ratio 33.]{
      \includegraphics[width=0.30\textwidth,clip]{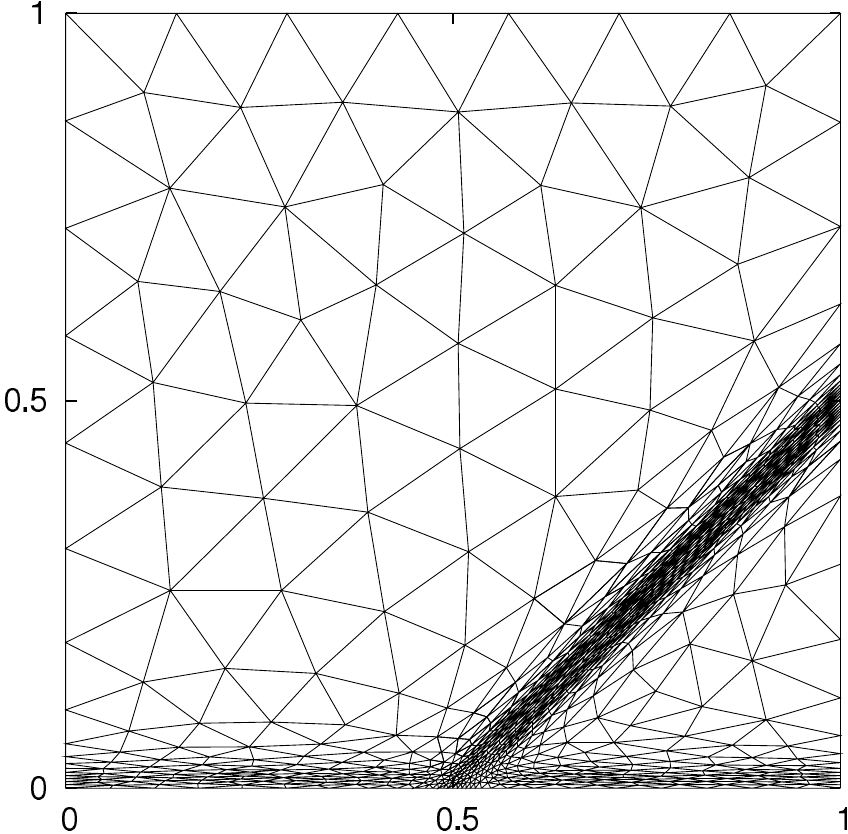}
      \qquad
      \includegraphics[width=0.31\textwidth,clip]{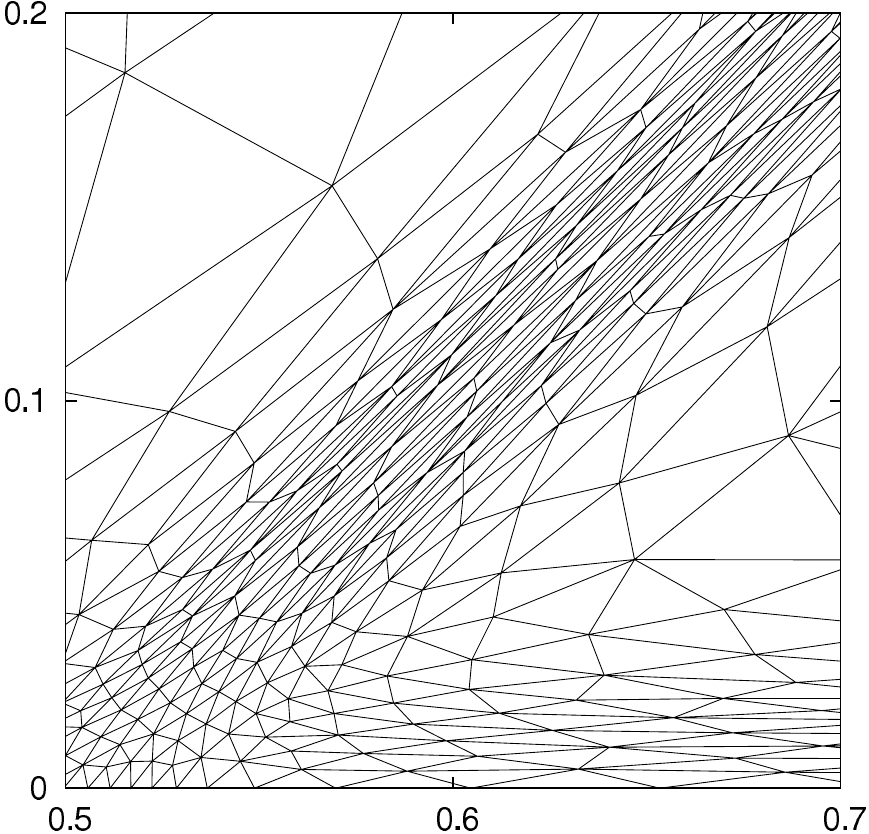}
      \label{fig:tanh:vfhbee}
   } \\
  \subfloat[HBEE: 677 vertices and 1256 triangles,
         $\norm{\nabla e_h}_{L^2} = 6.0 \times 10^{-1}$,
         maximum aspect ratio 45.]{
      \includegraphics[width=0.30\textwidth,clip]{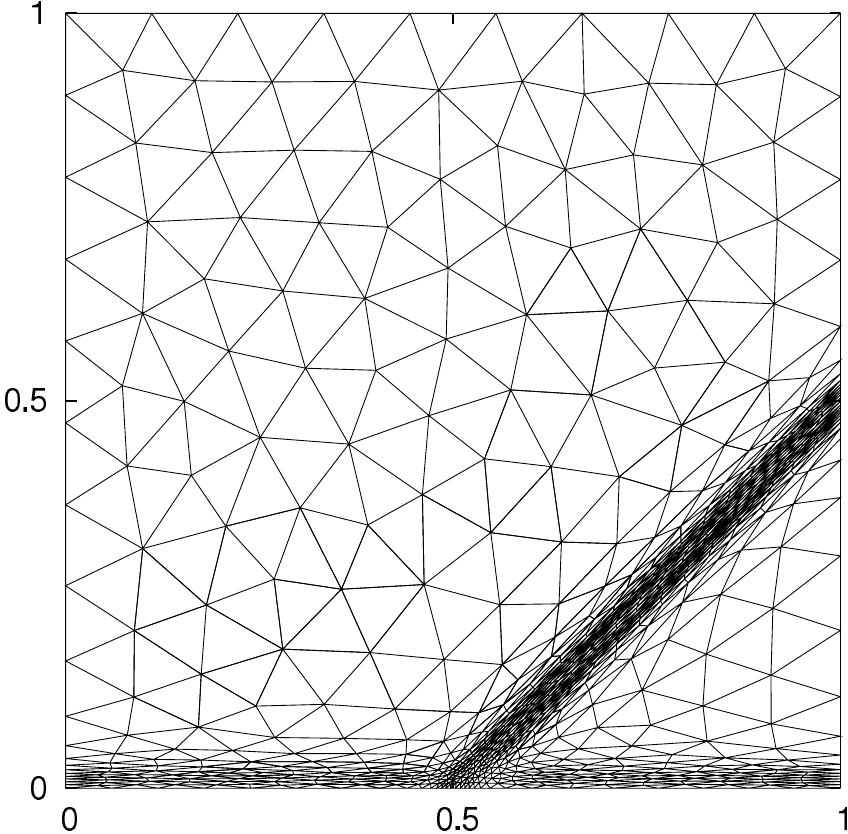}
      \qquad
      \includegraphics[width=0.31\textwidth,clip]{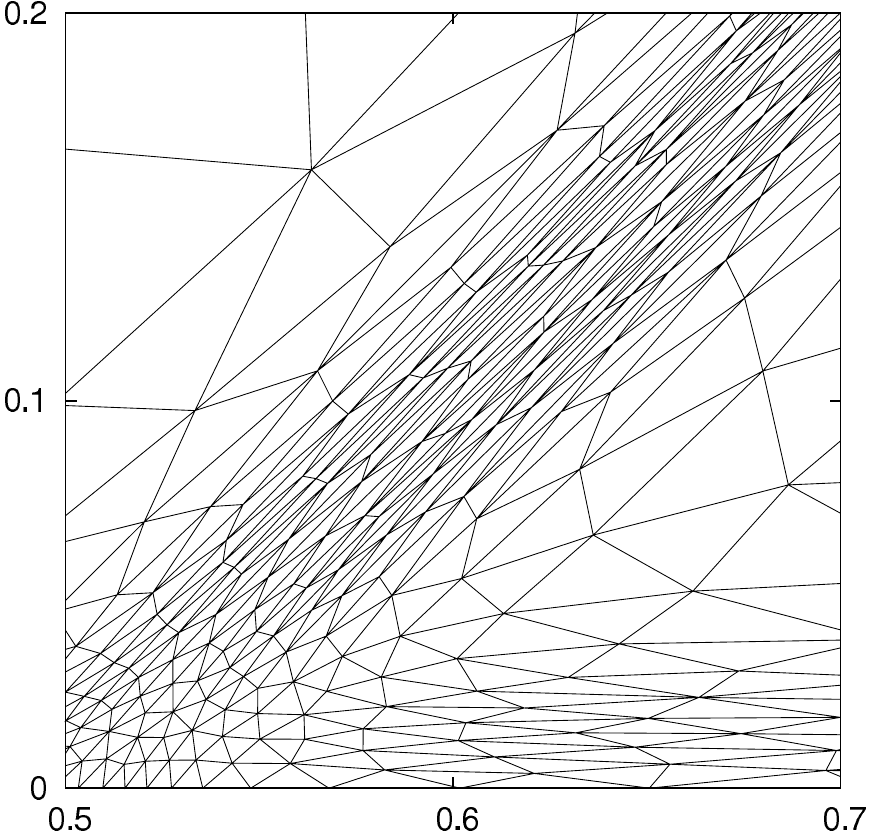}
      \label{fig:tanh:hbee}
   }
   \caption{Example~\ref{ex:tanh} (anisotropic):
      adaptive meshes obtained with different anisotropic metric tensors (left);
      close-up views at (0.5,0)(right).}
     \label{fig:tanh}
\end{figure}

Figure~\ref{fig:tanh:convergence} shows $\norm{\nabla e_h}_{L^2(\Omega)}$ (the finite element solution error measured in the $H^1$ semi-norm) against the number of elements for all metric tensors; adaptive mesh examples as well as close-up views near the point $(0.5,0)$, where the boundary layer meets the shock wave, are given in Figs.~\ref{fig:tanh:iso} and \ref{fig:tanh}.

% tanh example: qls<->hb: mesh examples
% beta = 0.8; nbv = 530 (qls 9), 530 (hb 9), 510 (hb 6)

First, we observe that all methods provide good mesh density, but anisotropic adaptation (Fig.~\ref{fig:tanh}) provides much better mesh alignment with the boundary layer and the shock wave front than isotropic adaptation (Fig.~\ref{fig:tanh:iso}).
Once again, it shows that anisotropic adaptation is advantageous and should be favored for such problems.
Meshes and the solution errors obtained with anisotropic adaptation are very similar.
Moreover, mesh adaptation by means of the variational formulation and HBEE is comparable with those obtained with both the Hessian recovery-based method in \cite{HuaLi10} and the PDE HBEE method in \cite{HuKaLa10}.

One may observe that the PDE HBEE method produces a slightly smaller error.
The reason for this behaviour can be attributed to the use of the residual and the edge jump in the metric tensor based on the variational formulation: it can be seen as a data smoothing, which prevents the sharp change of the element aspect ratio. 
The effect can be observed in Fig.~\ref{fig:tanh}: the maximum aspect ratio for the metric tensor based on the residual and HBEE is 33 (Fig.~\ref{fig:tanh:vfhbee}), whereas the maximum aspect ratio of the mesh obtained by means of just HBEE is 45 (Fig.~~\ref{fig:tanh:hbee}).
Hence,  the PDE HBEE based metric tensor might be a better choice for quadratic functionals with strong anisotropic features.

\subsection{A non-quadratic functional}
\label{sec:nonLinearExample}
\begin{figure}[t] \centering
   \includegraphics[width=0.38\textwidth,clip]{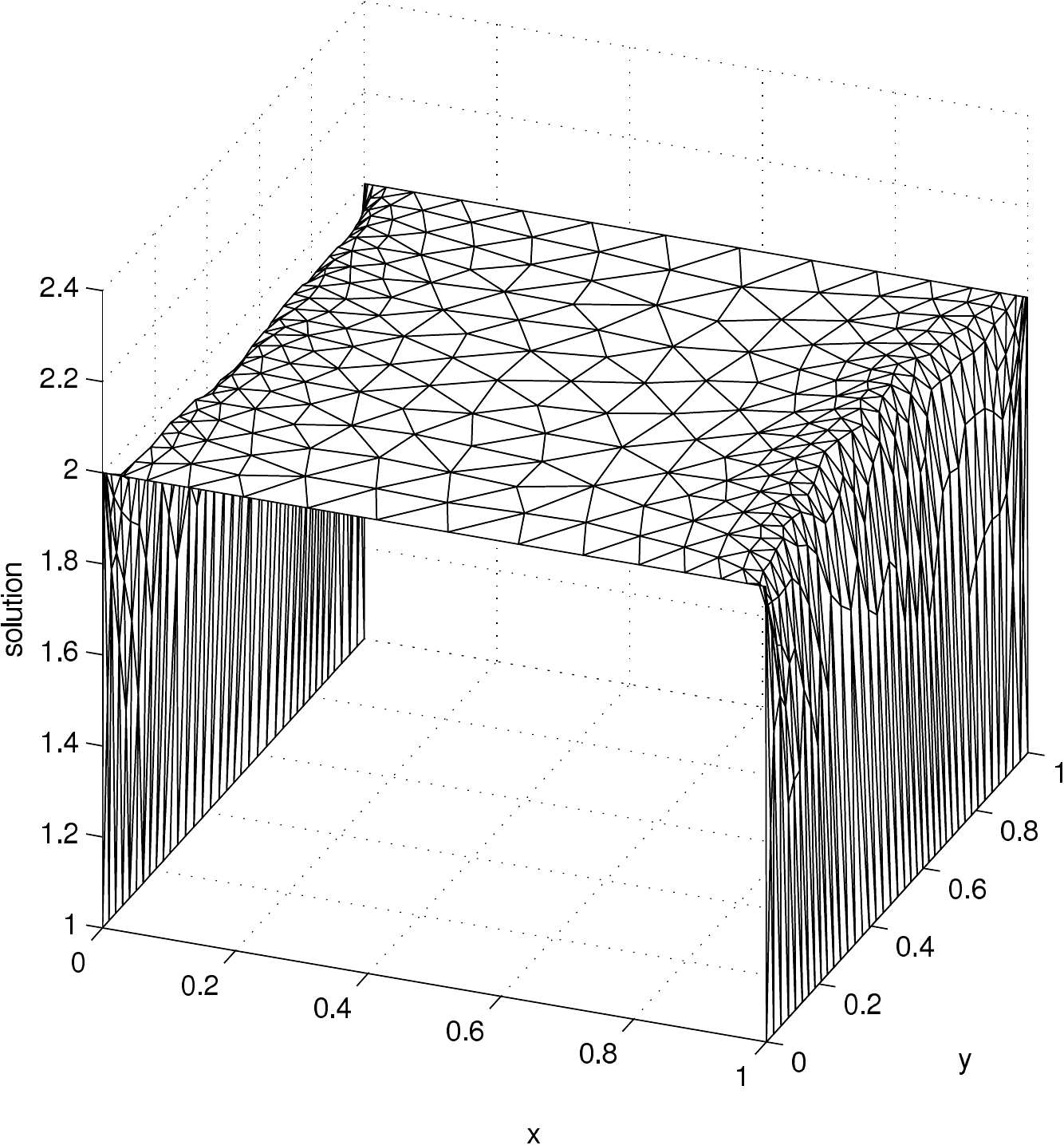}
   \caption{Example~\ref{sec:nonLinearExample}:
      surface plots of the numerical solutions.}
   \label{fig:plot:2}
\end{figure}

%--------------------------------------------------------------------------
The second example is an anisotropic variational problem \cite{Bildhauer03} defined by the functional
\begin{align}
   I[u] = \int_\Omega \left[ \left( 1 + \abs{\nabla u}^2 \right)^{3/4}
      + 1000 u_y^2 \right] d\boldx
   \label{ex3-fn}
\end{align}
with $\Omega = (0,1) \times (0,1)$ and the boundary condition 
\begin{align*}
   \begin{cases}
      u = 1  & \text{on } x=0 \text{ or } x=1, \\
      u = 2  & \text{on } y=0 \text{ or } y=1.
   \end{cases}
\end{align*}
Unlike the previous example, the functional \eqref{ex3-fn} is not quadratic.
Hence, the quantities $\abs{\delta I[u_h,e_h]}$ and $\left(\abs{\delta I[u_h,e_h]} / \abs{\nabla e_h}_{L^2(\Omega)}\right)^2$ are not mathematically equivalent to $\norm{\nabla e_h}^2_{L^2(\Omega )}$; this example is a good test for the formulas of the metric tensor based on the variational formulation.
Analytical solution is not available for this example.
A computed solution in Fig.~\ref{fig:plot:2} shows that the mesh adaptation challenge for this example lies  in the resolution of the sharp boundary layers near $x = 0$ and $x = 1$.

\begin{figure}[p] \centering
  \subfloat[Variational formulation with residual and edge jumps (isotropic):
         644 vertices and 1097 triangles, maximum aspect ratio 3.4.]{
      \includegraphics[width=0.30\textwidth,clip]{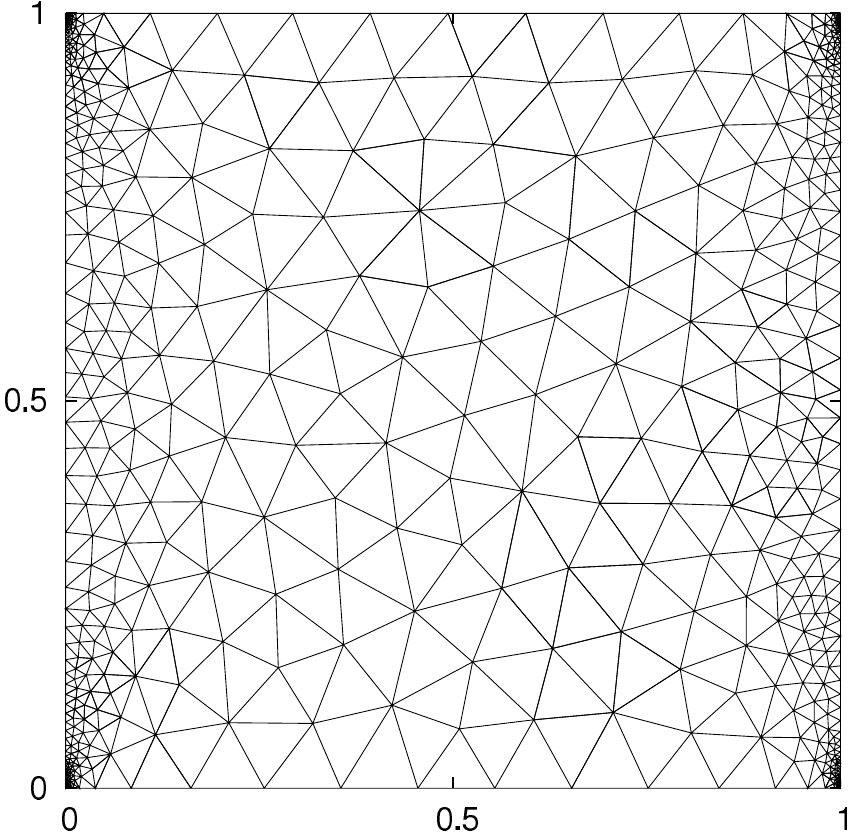}
      \qquad
      \includegraphics[width=0.31\textwidth,clip]{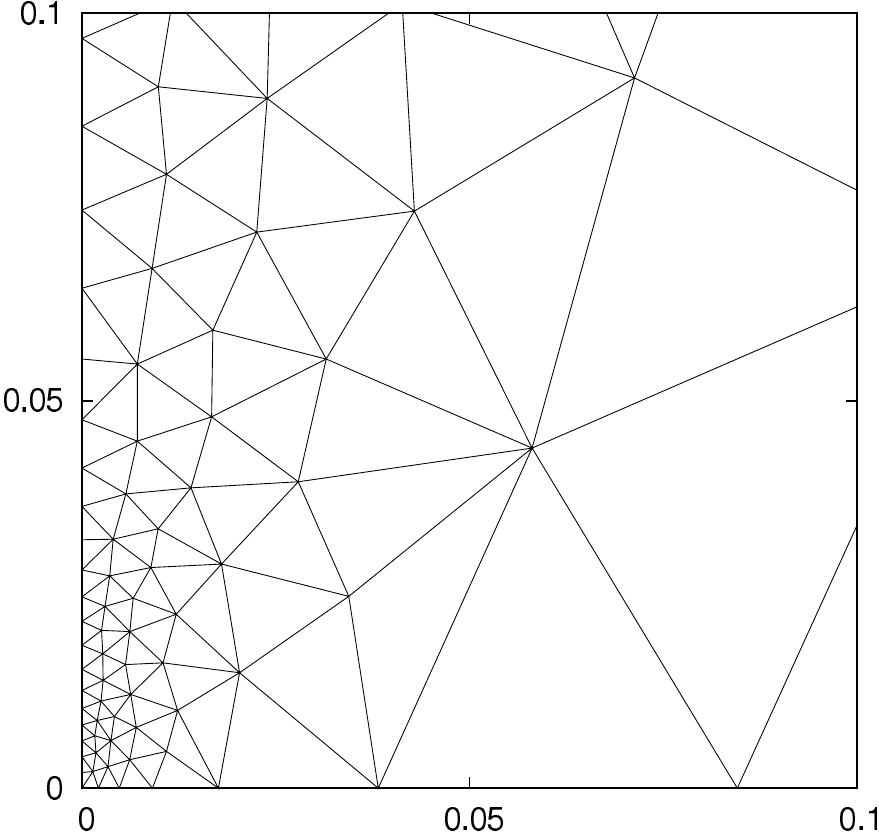}
      \label{fig:nonLinearExample:iso}
   }\\
   \subfloat[Variational formulation with residual, edge jumps, and Hessian recovery (anisotropic):
         649 vertices and 1160 triangles, maximum aspect ratio 15.]{
      \includegraphics[width=0.30\textwidth,clip]{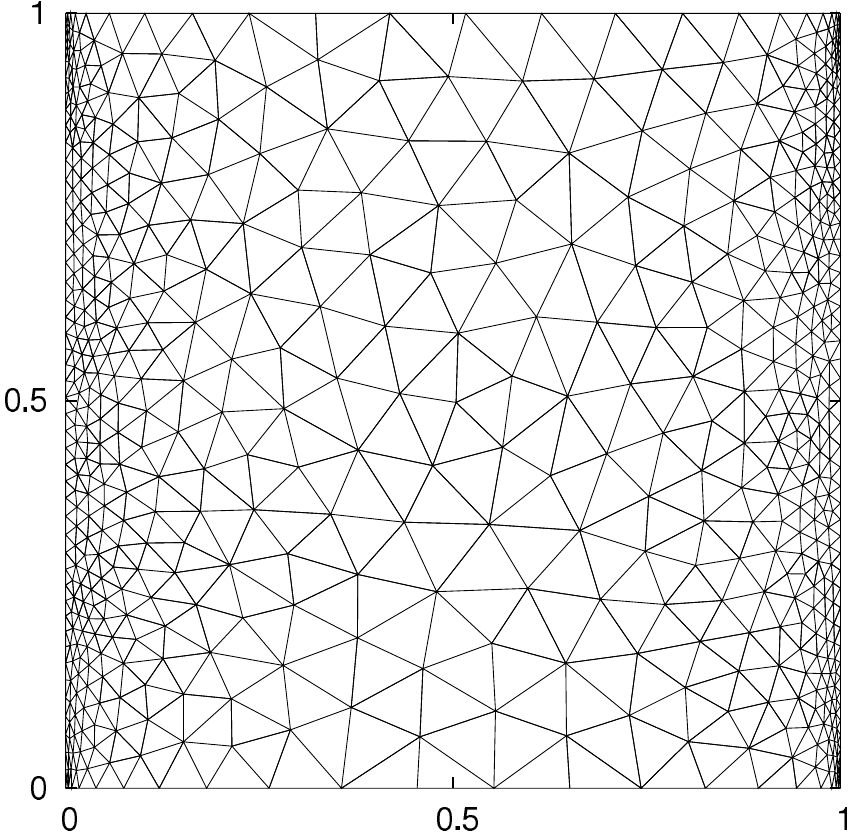}
      \qquad
      \includegraphics[width=0.31\textwidth,clip]{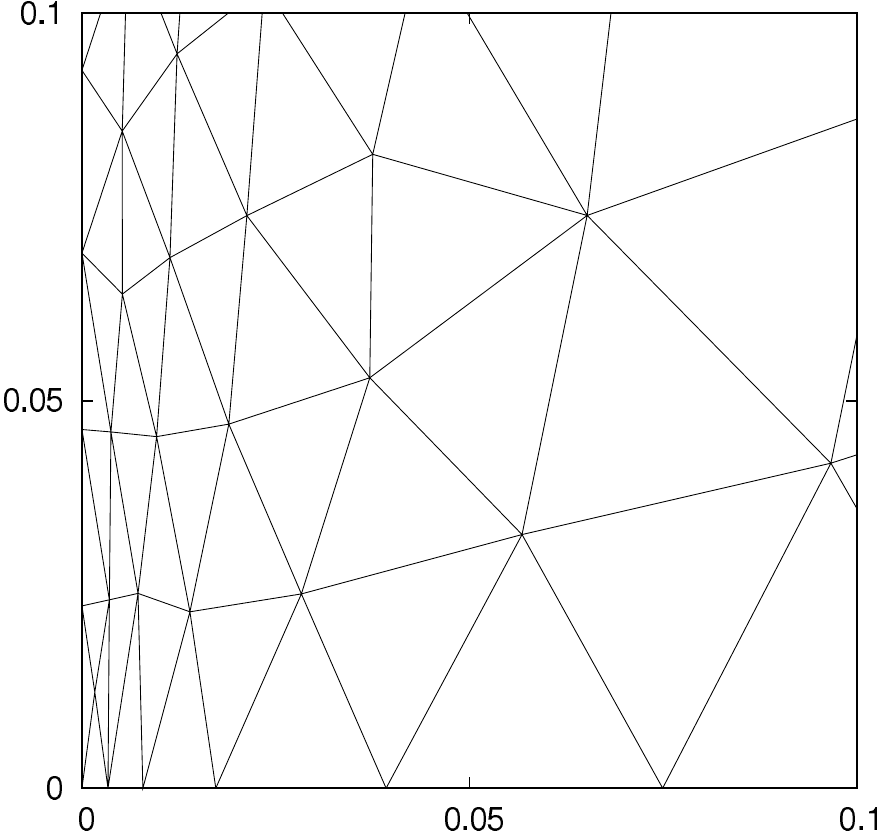}
      \label{fig:nonLinearExample:qls}
   }\\
   \subfloat[Variational formulation with residual, edge jumps, and HBEE:
            639 vertices and 1143 triangles, maximum aspect ratio 51.]{
      \includegraphics[width=0.30\textwidth,clip]{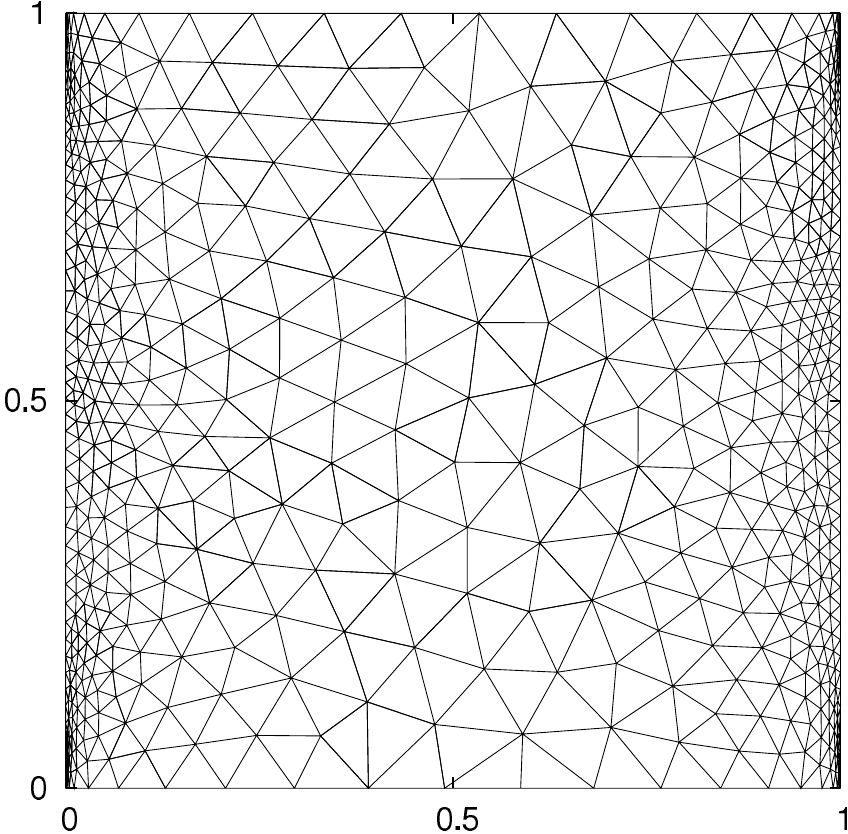}
      \qquad
      \includegraphics[width=0.31\textwidth,clip]{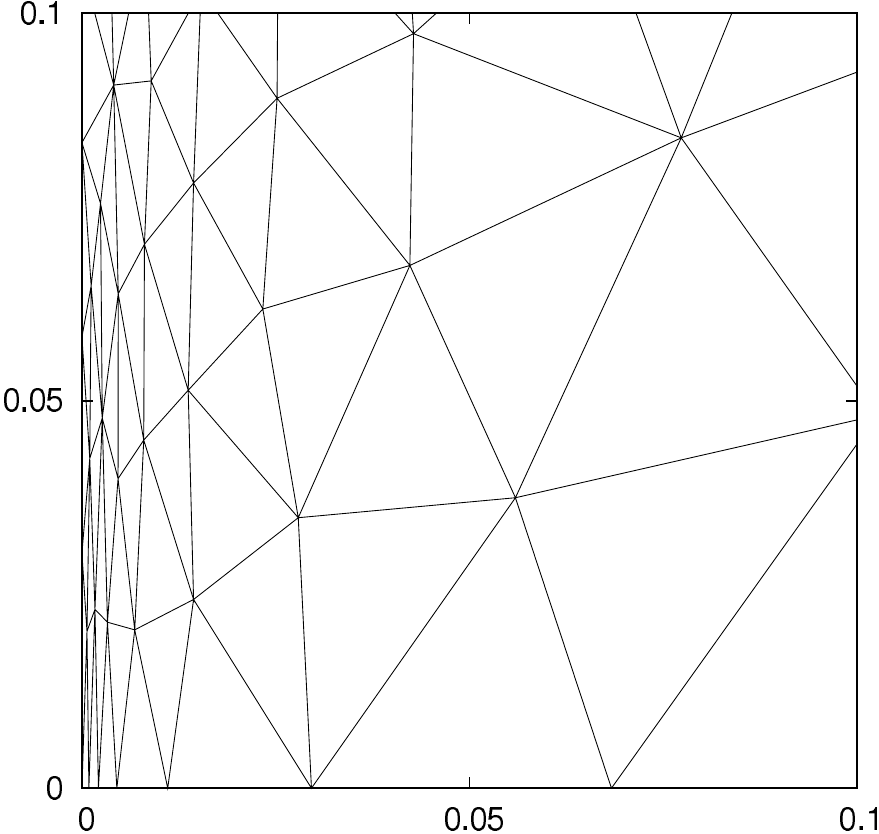}
      \label{fig:nonLinearExample:hb}
   }
     \caption{Example~\ref{sec:nonLinearExample}: adaptive meshes obtained 
        by means of isotropic, anisotropic Hessian recovery-based, 
        and anisotropic HBEE-based metric tensors (left); 
        close-up views at (0,0)(right).}
     \label{fig:nonLinearExample}
\end{figure}

Adaptive meshes obtained with the three different metric tensors (isotropic, Hessian recovery-based, and HBEE-based) are given in Fig.~\ref{fig:nonLinearExample}.
They all have correct mesh concentration, but the anisotropic metric tensors (Figs.~\ref{fig:nonLinearExample:qls} and \ref{fig:nonLinearExample:hb}) provide a much better alignment with the boundary layers.
Again, both anisotropic meshes are comparable, although mesh elements near the boundary layer in the HBEE-based adaptive mesh have a larger aspect ratio than elements of the mesh obtained by means of the Hessian recovery.
This could be due to the smoothing nature of the Hessian recovery: usually, it operates on a larger patch, thus introducing an additional smoothing effect, which affects the grading of the elements' size and orientation.
The global hierarchical basis error estimator does not have this handicap and, in this example, the mesh obtained by means of HBEE is slightly better aligned with the steep boundary layers.

% non-linear functional example: qls<->hb: mesh examples
% beta = 0.8; nbv = 650 (iso), 510 (qls), 500 (hb)

\subsection{Image processing}
\label{sec:imageProcessingExample}
\begin{figure}[t] \centering
   \includegraphics[width=0.38\textwidth,clip]{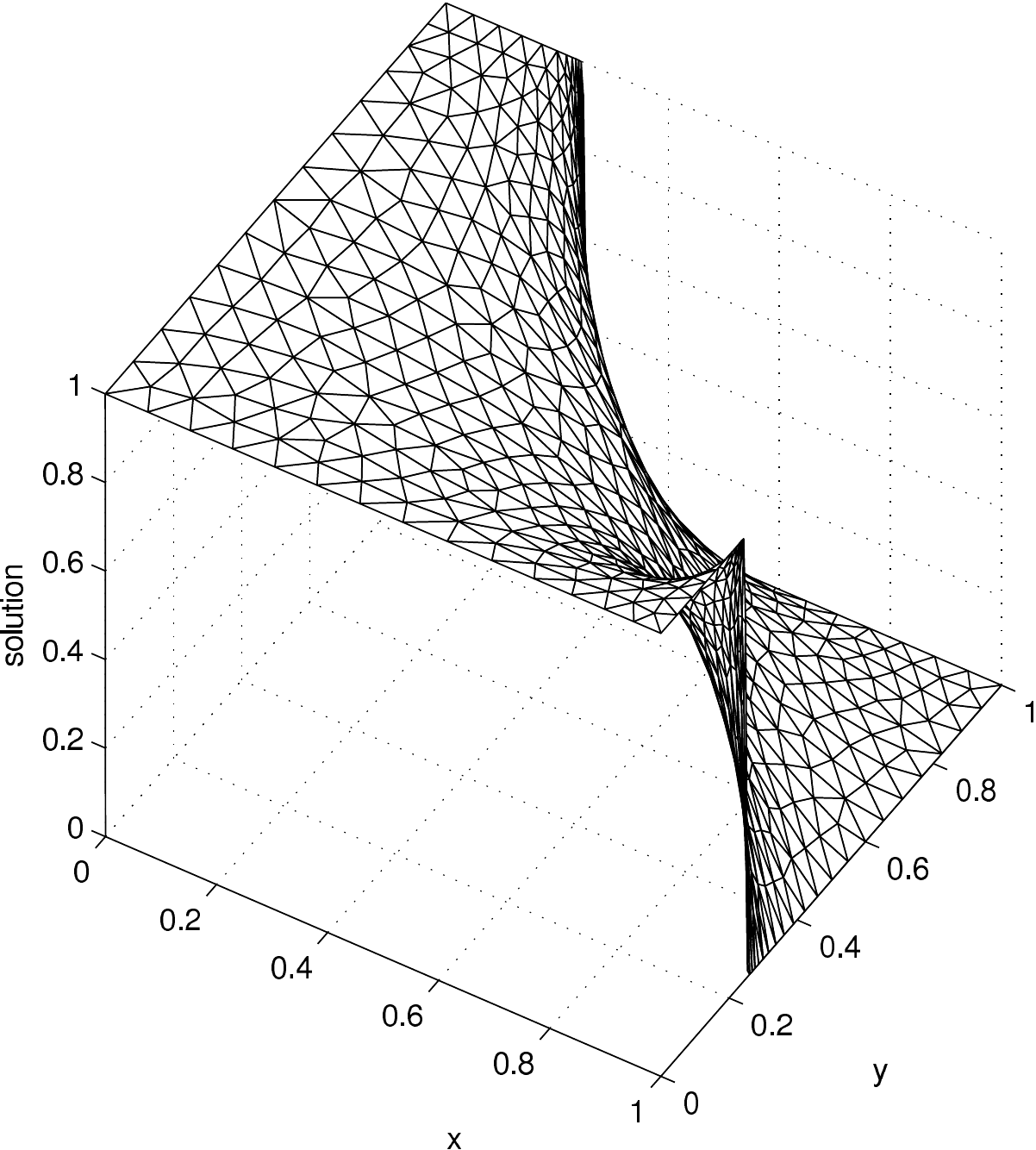}
   \caption{Example~\ref{sec:imageProcessingExample}:
      surface plots of the numerical solution.}
    \label{fig:plot:3}
\end{figure}

\begin{figure}[p] \centering
  \subfloat[Variational formulation with residual and edge jumps (isotropic):
            657 vertices and 1207 triangles, maximum aspect ratio 2.7.]{
      \includegraphics[width=0.30\textwidth,clip]{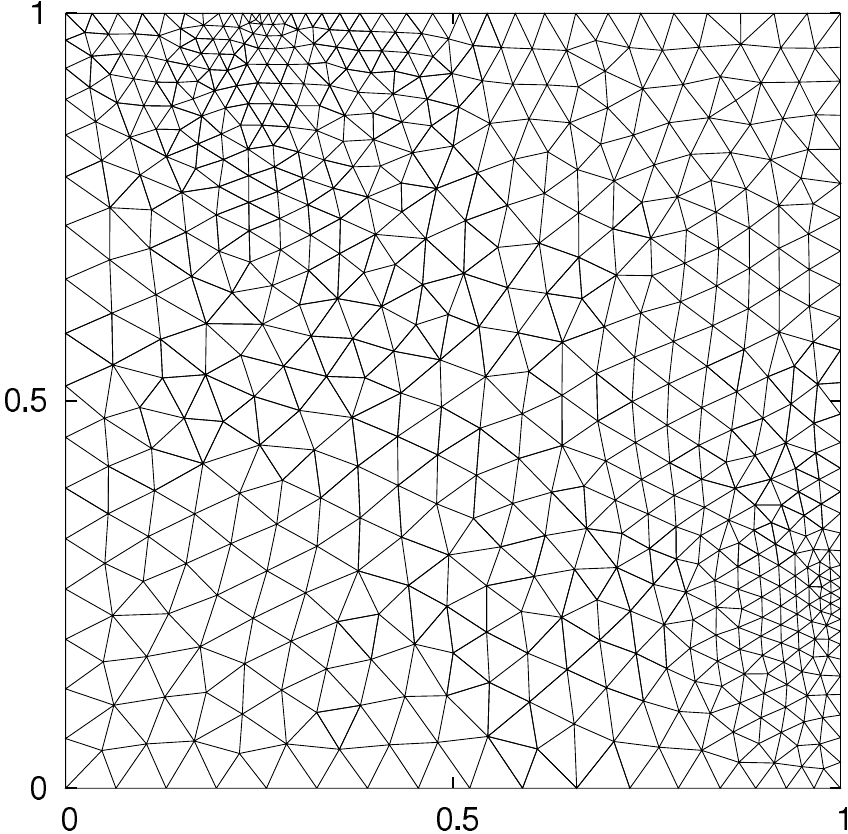}
      \qquad
      \includegraphics[width=0.31\textwidth,clip]{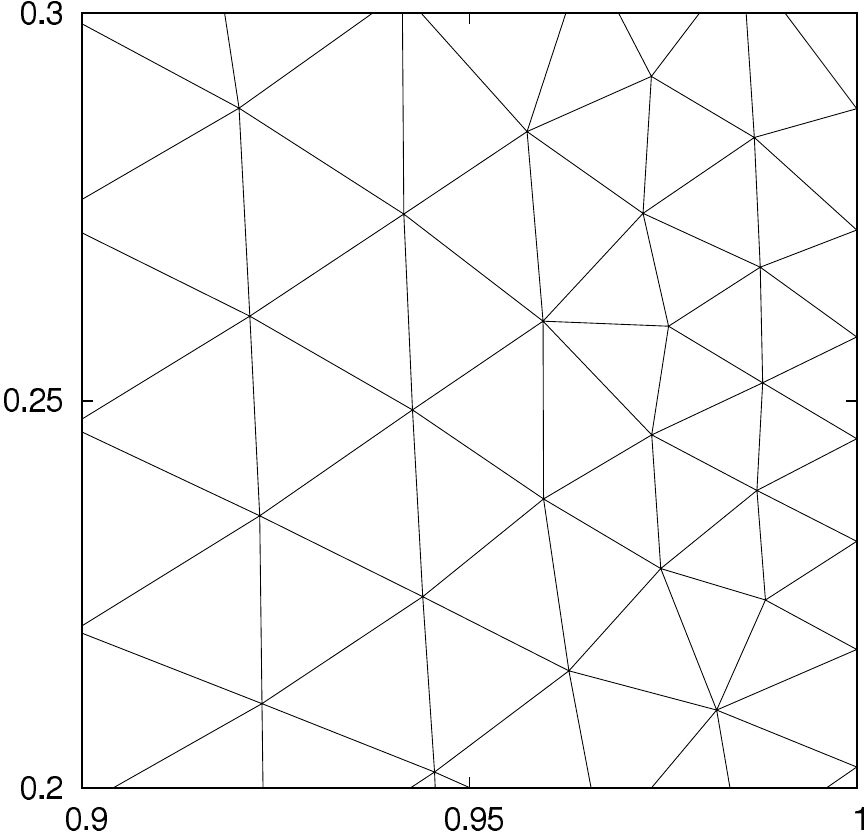}
   }\\
   \subfloat[Variational formulation with residual, edge jumps, and Hessian recovery (anisotropic):
            662 vertices and 1219 triangles, maximum aspect ratio 3.5.]{
      \includegraphics[width=0.30\textwidth,clip]{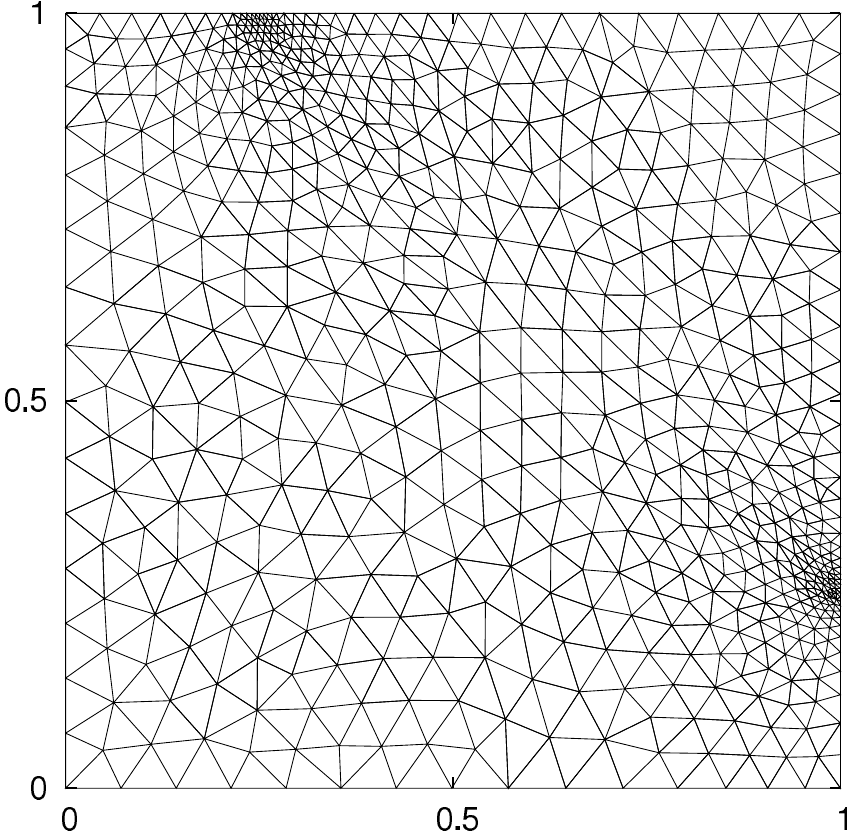}
      \qquad
      \includegraphics[width=0.31\textwidth,clip]{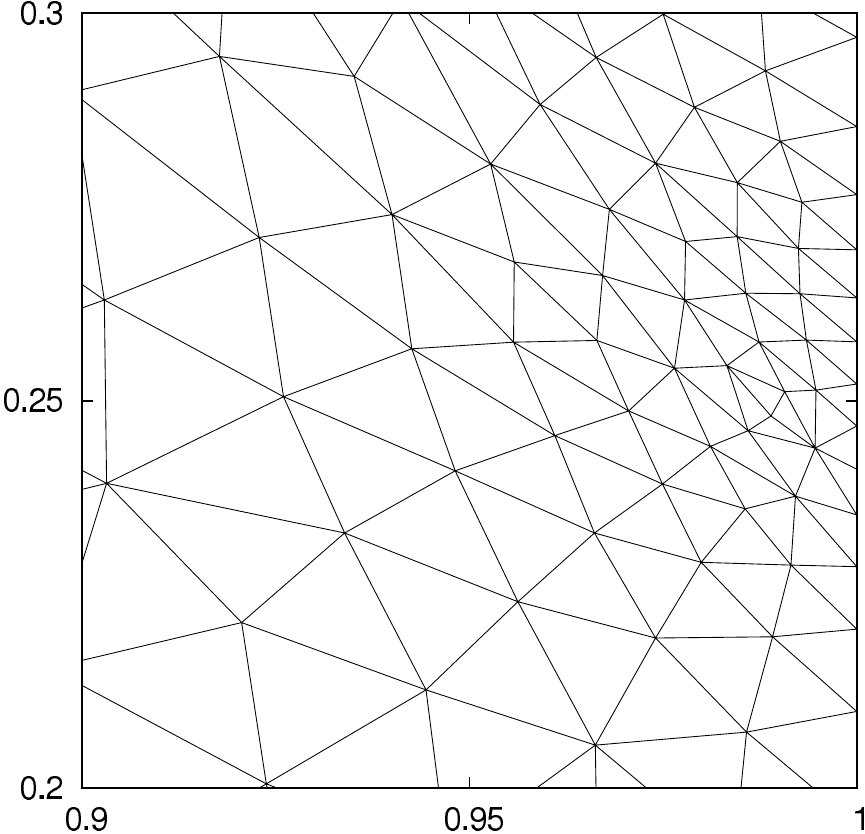}
   }\\
   \subfloat[Variational formulation with residual, edge jumps, and HBEE:
            656 vertices and 1209 triangles, maximum aspect ratio 4.7.]{
      \includegraphics[width=0.30\textwidth,clip]{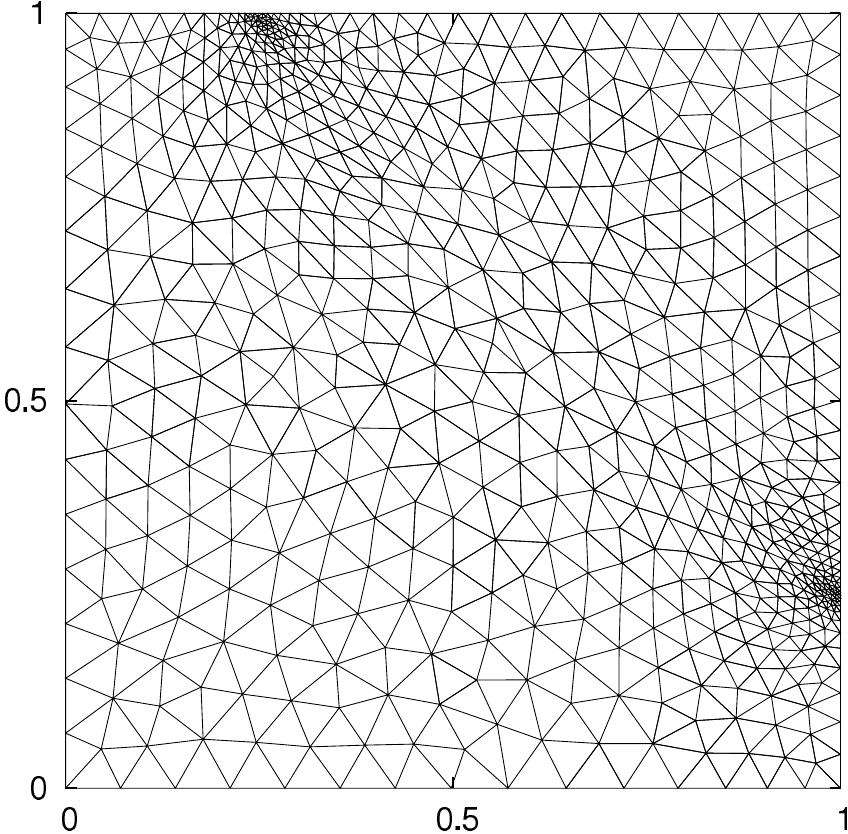}
      \qquad
      \includegraphics[width=0.31\textwidth,clip]{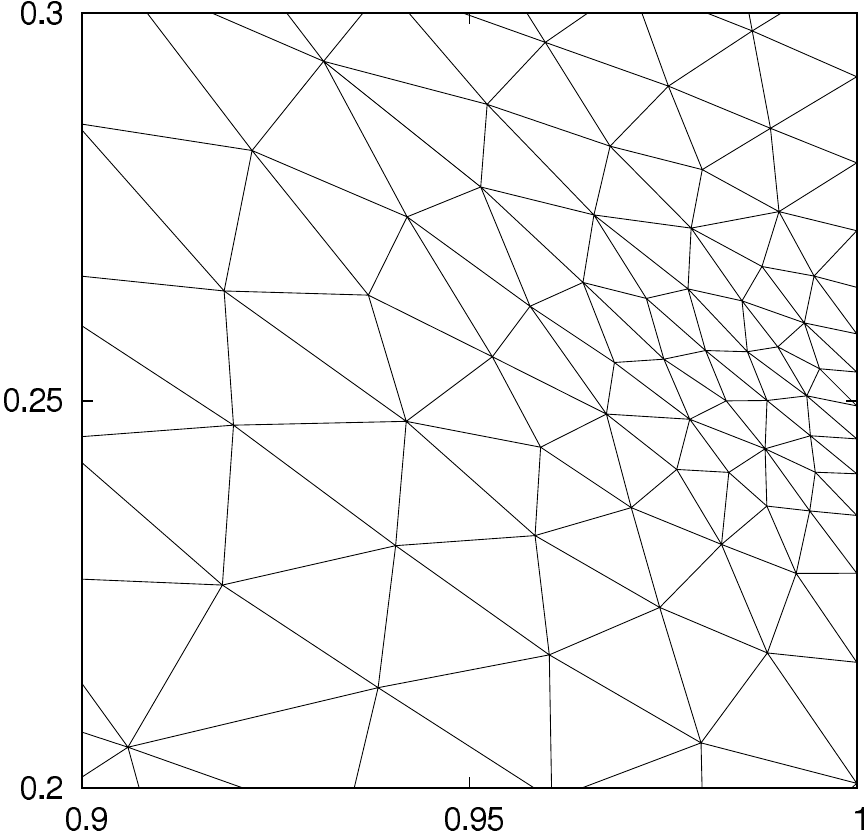}
   }
     \caption{Example~\ref{sec:imageProcessingExample}:  adaptive meshes 
        obtained by means of isotropic, anisotropic Hessian recovery-based,
        and anisotropic HBEE-based metric tensors (left);
        close-up views at (1,0.25)(right).}
     \label{fig:imageProcessing}
\end{figure}

%--------------------------------------------------------------------------
Our third example is an energy functional
\begin{align}
   I[u] = \int_\Omega \left[ \left(p - u\right)^2 
      + \left(1 + \abs{\nabla u}^2 \right)^{1/2} \right]  d\boldx,
   \label{ex4-fn}
\end{align}
which is used in image processing with the observed image $p$ and the reconstructed image $u$ \cite{AubVes97}. 
In our computation, we choose $\Omega = (0,1) \times (0,1)$,
\[ 
   p = 1 / \left(1 + e^{1000 (x + y - 1.25)}\right)
\]
and the boundary condition
\[
   u = p \quad \text{on } \partial \Omega.
\]
The analytical solution is not known; a computed numerical solution is shown in Fig.~\ref{fig:plot:3} and Fig.~\ref{fig:imageProcessing} shows examples of adaptive meshes.

% image processing example: qls<->hb: mesh examples
% beta = 0.8; nbv =560 (iso), 510 (qls), 495 (hb)

As in the previous example, the anisotropic metric tensors are comparable and provide a better mesh adaptation than the isotropic one.
Again, the HBEE-based mesh has a slightly larger maximum aspect ratio.

%**************************************************************************
\section{Concluding remarks}
%**************************************************************************
\label{sect:Conclusion}

Numerical results confirm the conclusion of \cite{HuKaLa10} that a global HBEE can be a successful alternative to Hessian recovery in mesh adaptation; a fast approximate solution of the global error problem is sufficient to provide directional information for anisotropic mesh adaptation. 
They also confirm the conjecture that good mesh adaptation does not require a convergent Hessian recovery or an accurate error estimator, but rather some additional information of global nature, although it is still unclear which information exactly is necessary for successful anisotropic adaptation.

For quadratic functionals, numerical results suggest that algorithms which involve non-directional information such as residual sometime can be slightly disadvantageous for problems with sharp anisotropic features.
Such algorithms can smooth out the metric tensor and therefore prevent the rapid change of element size and orientation, which is highly desired for this class of problems.
Also, as mentioned in \cite{HuKaLa10}, Hessian recovery could result in unnecessarily high mesh density for problems with discontinuities.
Global hierarchical error estimation seems less affected by these issues and, depending on the underlying problem, can provide a more robust solution.

\vspace{1em}

\textbf{Acknowledgment.} This work was supported in part by the National Science Foundation (U.S.A.) under grant DMS-0712935 and by the German Research Foundation (DFG) under grant KA~3215/1-1.
The authors are grateful to the anonymous referees for their valuable comments.

%\vspace{2em}
\newpage
%**************************************************************************
%\bibliographystyle{model1b-num-names}
%\bibliography{HuKaLi10}

\begin{thebibliography}{16}
\expandafter\ifx\csname natexlab\endcsname\relax\def\natexlab#1{#1}\fi
\providecommand{\bibinfo}[2]{#2}
\ifx\xfnm\relax \def\xfnm[#1]{\unskip,\space#1}\fi
%Type = Article
\bibitem[{Aubert and Vese(1997)}]{AubVes97}
\bibinfo{author}{G.~Aubert}, \bibinfo{author}{L.~Vese}, \bibinfo{title}{A
  variational method in image recovery}, \bibinfo{journal}{SIAM Journal on
  Numerical Analysis} \bibinfo{volume}{34} (\bibinfo{year}{1997})
  \bibinfo{pages}{1948--1979}.
%Type = Article
\bibitem[{Bank and Smith(1993)}]{BanSmi93}
\bibinfo{author}{R.E. Bank}, \bibinfo{author}{R.K. Smith}, \bibinfo{title}{A
  posteriori error estimates based on hierarchical bases},
  \bibinfo{journal}{SIAM J. Numer. Anal.} \bibinfo{volume}{30}
  (\bibinfo{year}{1993}) \bibinfo{pages}{921--935}.
%Type = Article
\bibitem[{Becker and Rannacher(1996)}]{BR96}
\bibinfo{author}{R.~Becker}, \bibinfo{author}{R.~Rannacher}, \bibinfo{title}{A
  feed-back approach to error control in finite element methods: basic analysis
  and examples}, \bibinfo{journal}{East-West J. Numer. Math.}
  \bibinfo{volume}{4} (\bibinfo{year}{1996}) \bibinfo{pages}{237--264}.
%Type = Article
\bibitem[{Becker and Rannacher(2001)}]{BR01}
\bibinfo{author}{R.~Becker}, \bibinfo{author}{R.~Rannacher}, \bibinfo{title}{An
  optimal control approach to a posteriori error estimation in finite element
  methods}, \bibinfo{journal}{Acta Numer.} \bibinfo{volume}{10}
  (\bibinfo{year}{2001}) \bibinfo{pages}{1--102}.
%Type = Article
\bibitem[{Berndt et~al.(1997)Berndt, Manteuffel and Mccormick}]{BMM97}
\bibinfo{author}{M.~Berndt}, \bibinfo{author}{T.A. Manteuffel},
  \bibinfo{author}{S.F. Mccormick}, \bibinfo{title}{Local error estimates and
  adaptive refinement for first-order system least squares ({FOSLS})},
  \bibinfo{journal}{Electr. Trans. Numer. Anal.} \bibinfo{volume}{6}
  (\bibinfo{year}{1997}) \bibinfo{pages}{35--43}.
%Type = Book
\bibitem[{Bildhauer(2003)}]{Bildhauer03}
\bibinfo{author}{M.~Bildhauer}, \bibinfo{title}{Convex variational problems.
  Linear, nearly linear and anisotropic growth conditions.}, volume
  \bibinfo{volume}{1818} of \textit{\bibinfo{series}{Lecture Notes in
  Mathematics}}, \bibinfo{publisher}{Springer Berlin / Heidelberg},
  \bibinfo{year}{2003}.
%Type = Article
\bibitem[{Borouchaki et~al.(1997{\natexlab{a}})Borouchaki, George, Hecht, Laug
  and Saltel}]{BoGHLS97}
\bibinfo{author}{H.~Borouchaki}, \bibinfo{author}{P.L. George},
  \bibinfo{author}{F.~Hecht}, \bibinfo{author}{P.~Laug},
  \bibinfo{author}{E.~Saltel}, \bibinfo{title}{Delaunay mesh generation
  governed by metric specifications. {P}art {I}. {A}lgorithms},
  \bibinfo{journal}{Finite Elements in Analysis and Design}
  \bibinfo{volume}{25} (\bibinfo{year}{1997}{\natexlab{a}})
  \bibinfo{pages}{61--83}.
%Type = Article
\bibitem[{Borouchaki et~al.(1997{\natexlab{b}})Borouchaki, George and
  Mohammadi}]{BoGeMo97}
\bibinfo{author}{H.~Borouchaki}, \bibinfo{author}{P.L. George},
  \bibinfo{author}{B.~Mohammadi}, \bibinfo{title}{Delaunay mesh generation
  governed by metric specifications. {P}art {II}. {A}pplications},
  \bibinfo{journal}{Finite Elem. Anal. Des.} \bibinfo{volume}{25}
  (\bibinfo{year}{1997}{\natexlab{b}}) \bibinfo{pages}{85--109}.
%Type = Article
\bibitem[{Castro-D{\'\i}az et~al.(1997)Castro-D{\'\i}az, Hecht, Mohammadi and
  Pironneau}]{CaHeMP97}
\bibinfo{author}{M.J. Castro-D{\'\i}az}, \bibinfo{author}{F.~Hecht},
  \bibinfo{author}{B.~Mohammadi}, \bibinfo{author}{O.~Pironneau},
  \bibinfo{title}{Anisotropic unstructured mesh adaption for flow simulations},
  \bibinfo{journal}{Int. J. Numer. Meth. Fluids} \bibinfo{volume}{25}
  (\bibinfo{year}{1997}) \bibinfo{pages}{475--491}.
%Type = Article
\bibitem[{Deuflhard et~al.(1989)Deuflhard, Leinen and Yserentant}]{DeLeYs89}
\bibinfo{author}{P.~Deuflhard}, \bibinfo{author}{P.~Leinen},
  \bibinfo{author}{H.~Yserentant}, \bibinfo{title}{Concepts of an adaptive
  hierarchical finite element code}, \bibinfo{journal}{Impact Comput. Sci.
  Engrg.} \bibinfo{volume}{1} (\bibinfo{year}{1989}) \bibinfo{pages}{3--35}.
%Type = Article
\bibitem[{D{\"o}rfler and Nochetto(2002)}]{DoeNoc02}
  \bibinfo{author}{W. D{\"o}rfler}, \bibinfo{author}{R.H. Nochetto},
  \bibinfo{title}{Small data oscillation implies the saturation assumption},
  \bibinfo{journal}{Numer. Math.}
  \bibinfo{volume}{91} (\bibinfo{year}{2002}) \bibinfo{pages}{1--12}.
%Type = Book
\bibitem[{Evans(1998)}]{Evans98}
\bibinfo{author}{L.C. Evans}, \bibinfo{title}{Partial Differential Equations},
  volume~\bibinfo{volume}{19} of \textit{\bibinfo{series}{Graduate Studies in
  Mathematics}}, \bibinfo{publisher}{American Mathematical Society, Providence,
  Rhode Island}, \bibinfo{year}{1998}.
%Type = Book
\bibitem[{Ern and Guermond(2004)}]{Ern04}
  \bibinfo{author}{A. Ern}, \bibinfo{author}{J.-L. Guermond}, 
  \bibinfo{title}{Theory and Practice of Finite Elements},
  volume~\bibinfo{volume}{159} of \textit{\bibinfo{series}{Applied Mathematical
  Sciences}}, \bibinfo{publisher}{Springer-Verlag New York}, \bibinfo{year}{2004}.
%Type = Book
\bibitem[{Frey and George(2008)}]{Frey08}
\bibinfo{author}{P.J. Frey}, \bibinfo{author}{P.L. George},
  \bibinfo{title}{Mesh Generation. Second Edition}, \bibinfo{publisher}{John
  Wiley \& Sons, Inc., Hoboken, NJ}, \bibinfo{year}{2008}.
%Type = Article
\bibitem[{Huang(2005)}]{Huang05a}
\bibinfo{author}{W.~Huang}, \bibinfo{title}{Metric tensors for anisotropic mesh
  generation}, \bibinfo{journal}{J. Comput. Phys.} \bibinfo{volume}{204}
  (\bibinfo{year}{2005}) \bibinfo{pages}{633--665}.
%Type = Article
\bibitem[{Huang(2006)}]{Huang06}
\bibinfo{author}{W.~Huang}, \bibinfo{title}{Mathematical principles of
  anisotropic mesh adaptation}, \bibinfo{journal}{Commun. Comput. Phys.}
  \bibinfo{volume}{1} (\bibinfo{year}{2006}) \bibinfo{pages}{276--310}.
%Type = Article
\bibitem[{Huang et~al.(2010)Huang, Kamenski and Lang}]{HuKaLa10}
\bibinfo{author}{W.~Huang}, \bibinfo{author}{L.~Kamenski},
  \bibinfo{author}{J.~Lang}, \bibinfo{title}{A new anisotropic mesh adaptation
  method based upon hierarchical a posteriori error estimates},
  \bibinfo{journal}{J. Comput. Phys.} \bibinfo{volume}{229}
  (\bibinfo{year}{2010}) \bibinfo{pages}{2179--2198}.
%Type = Article
\bibitem[{Huang and Li(2010)}]{HuaLi10}
\bibinfo{author}{W.~Huang}, \bibinfo{author}{X.~Li}, \bibinfo{title}{An
  anisotropic mesh adaptation method for the finite element solution of
  variational problems}, \bibinfo{journal}{Finite Elem. Anal. Des.}
  \bibinfo{volume}{46} (\bibinfo{year}{2010}) \bibinfo{pages}{61--73}.
\end{thebibliography}

%**************************************************************************

\end{document}